\documentclass[11pt,a4paper]{article}

\usepackage[T1]{fontenc}
\usepackage[utf8]{inputenc}
\usepackage{lmodern,amsmath,amssymb,amsthm,mathrsfs,microtype}
\usepackage[a4paper,left=25mm,right=25mm,top=25mm,bottom=25mm]{geometry}
\usepackage{xcolor}
\usepackage[bookmarks=true,colorlinks=true,linkcolor=blue!35!black,citecolor=blue!35!black,urlcolor=blue!50!black]{hyperref}

\numberwithin{equation}{section}
\newtheorem{theorem}{Theorem}[section]
\newtheorem{proposition}[theorem]{Proposition}
\newtheorem{lemma}[theorem]{Lemma}
\newtheorem{corollary}[theorem]{Corollary}
\theoremstyle{remark}
\newtheorem{remark}[theorem]{Remark}

\hypersetup{pdftitle={The enhanced polylogarithm conjecture in weight two for imaginary quadratic fields},pdfauthor={Nobuo Sato},pdfsubject={2020 MSC: Primary 11R42; Secondary 11G15, 11G55, 11M41, 14C25, 19F27},pdfkeywords={Shintani zeta functions, enhanced dilogarithm, ray class fields, complex multiplication, higher Chow groups, regulators}}
\begin{document}
\begin{center}
{\LARGE The enhanced polylogarithm conjecture\\[3pt]in weight two for imaginary quadratic fields\par}
\vspace{12pt}
{\large Nobuo Sato\par}
\vspace{5pt}
{\small Department of Mathematics, National Taiwan University\\
Taipei 106319, Taiwan\\
\href{mailto:nbsato@ntu.edu.tw}{\texttt{nbsato@ntu.edu.tw}}\par}
\end{center}
\medskip
\begin{abstract}
Zagier's polylogarithm conjecture predicts that special values of Dedekind
zeta functions can be expressed in terms of polylogarithms at algebraic
arguments. For imaginary quadratic fields, Zagier and Gangl refined it by
lifting both sides to complex values modulo rational multiples of a power
of $2\pi i$. In previous work, the author extended their enhanced zeta
values from ideal classes to ray classes by means of the Shintani
L-function of two variables, and formulated the enhanced polylogarithm
conjecture for ray classes. In this article, we prove this conjecture in
weight two. More precisely, we establish a Galois-equivariant equality
between the enhanced partial zeta values of ray classes and the enhanced
dilogarithms of elements of the rational Bloch groups of the corresponding
ray class fields. These elements are constructed explicitly, starting from
Doran--Kerr's Eisenstein class on the square of a CM elliptic curve and
contracting it with the graphs of complex multiplication.
\end{abstract}



\section{Introduction}

The main theme of this article is a refinement of Zagier's polylogarithm
conjecture for imaginary quadratic fields. Zagier's conjecture relates
special values of Dedekind zeta functions to polylogarithms at algebraic
arguments \cite{Zagier1991}, and its dilogarithm case is naturally
expressed through the Bloch group and the regulator on algebraic $K_3$.

To explain the idea of the enhancement, let us begin with the classical
case of weight one. For a nonzero integral ideal $\mathfrak f\ne\mathcal O_K$
of an imaginary quadratic field $K$ and a ray class $A$ modulo $\mathfrak f$,
the Kronecker limit formula and complex multiplication give a Taylor
expansion of the partial zeta function of the form
\[
 \zeta_{\mathfrak f}(s,A)
   =\frac{2}{m}\log|\varepsilon_A|\,s+O(s^2),
\]
where $\varepsilon_A$ is a nonzero algebraic number in the ray class
field and $m$ is a positive integer
\cite[\S2, pp.~278--281]{Stark1977}.
The first Taylor coefficient records the absolute value of
$\varepsilon_A$, but loses its argument.

In weight one, the enhancement is a classical complex refinement of this
expansion. We describe it in terms of the Shintani L-function of two
variables considered in \cite{Sato}, since this description extends to
higher weights. We fix its fundamental domain and weighting, and write
$\mathscr L(s_1,s_2)$ for its normalized form, so that
\[
 \mathscr L(s,s)=\zeta_{\mathfrak f}(s,A).
\]
To illustrate the idea, suppose for simplicity that $\mathscr L$ is
holomorphic at the origin. For a suitable choice of $\varepsilon_A$ and
$m$, and of the branch of the logarithm, which depends on the fundamental
domain and weighting of $\mathscr L(s_1,s_2)$, its Taylor expansion takes
the form
\[
 \mathscr L(s_1,s_2)
   =\frac{\log\varepsilon_A}{m}\,s_1
    +\frac{\overline{\log\varepsilon_A}}{m}\,s_2
    +O\bigl((|s_1|+|s_2|)^2\bigr).
\]
On the diagonal, the two linear terms combine because
$\log\varepsilon_A+\overline{\log\varepsilon_A}=2\log|\varepsilon_A|$,
recovering the expansion above. In general, $\mathscr L$ need not be
holomorphic at the origin, but its partial derivatives have limits along
the diagonal, which take the place of the Taylor coefficients
\cite[Lemma~2.12 and Remark~2.26]{Sato}: differentiating the diagonal
identity for $s\ne0$ and taking the limit shows that their sum is still
$\zeta_{\mathfrak f}'(0,A)$. The refinement keeps the two limits, $(\log\varepsilon_A)/m$ and its
complex conjugate, separately rather than their sum. The Fourier
expressions for these limits \cite[Proposition~3.6]{Sato} are complex
logarithms of normalized theta functions. The modular
transformation laws of these functions, together with complex
multiplication, show that their values at the point of complex
multiplication are algebraic. This is how Stark obtains $\varepsilon_A$
\cite[\S2, pp.~278--281]{Stark1977}, and taking real parts recovers the
logarithmic formula above.

The situation in weight two is similar, with the complex logarithm
replaced by the enhanced dilogarithm. For the ordinary regulator
comparison, Deninger constructed $K$-theory classes associated with
Hecke characters of imaginary quadratic fields and computed their
regulators \cite{DeningerI,DeningerII}. Levin obtained explicit formulas
expressing Kronecker double series at complex multiplication points
as finite sums of Bloch--Wigner dilogarithms. At zero torsion argument,
these series give normalized partial zeta values at $s=2$ for ideal
classes of imaginary quadratic orders; his formulas also cover nonzero
torsion arguments
\cite{Levin}.

In 2000, Zagier and Gangl gave a refinement of Zagier's conjecture for
imaginary quadratic fields \cite[\S\S4 and 7]{GZ}. For each ideal class,
they constructed the enhanced zeta value, and on Bloch groups the
enhanced polylogarithm; both take values in $\mathbb C$ modulo rational
multiples of a power of $2\pi i$. They conjectured that the enhanced zeta
values are values of the enhanced polylogarithm at elements of the
rational Bloch group of the Hilbert class field.

In \cite{Sato}, the author extended the analytic side of Zagier--Gangl's
construction from ideal classes to ray classes, by means of the Shintani
L-function of two variables. For ideal classes, the resulting
invariants agree with Zagier--Gangl's enhanced zeta values, after
normalization and inversion of the ideal class \cite[Theorem~2]{Sato},
and \cite[Conjecture~4.1]{Sato} extends their conjecture to arbitrary ray
classes. The purpose of this article is to prove this conjecture in
weight two.

Fix an embedding of an imaginary quadratic field $K$ into $\mathbb C$.
For a nonzero integral ideal $\mathfrak f$, let $F=K_{\mathfrak f}$ be
the ray class field, and fix an embedding $F\hookrightarrow\mathbb C$
extending that of $K$.
The ray class group
$\operatorname{Cl}_{\mathfrak f}(K)$ is the group of fractional ideals
prime to $\mathfrak f$, modulo principal ideals $(c)$ with
$c\equiv1\pmod{\mathfrak f}$ at all primes dividing $\mathfrak f$.
Let
$\operatorname{rec}:\operatorname{Cl}_{\mathfrak f}(K)\to\operatorname{Gal}(F/K)$
be the reciprocity isomorphism. For a prime $\mathfrak p\nmid\mathfrak f$,
we normalize $\operatorname{rec}$ so that $\operatorname{rec}([\mathfrak p])^{-1}$
induces $x\mapsto x^{N(\mathfrak p)}$ on the residue fields of the primes of
$F$ above $\mathfrak p$.
This is the inverse of the usual Frobenius normalization, as in the
examples of \cite[\S4.3]{Sato}.

Put $\mathbb Q(j)=(2\pi i)^j\mathbb Q$.
For a ray class $A$, the enhanced partial zeta value
$\Lambda_1(-1,A)\in\mathbb C/\mathbb Q(1)$ is defined, as in weight one,
by the normalized diagonal limit at $(-1,-1)$ of the partial derivative
of the Shintani L-function with respect to $s_1$.
For $\sigma\in\operatorname{Gal}(F/K)$, we write
$\Lambda_1(-1,\sigma)=\Lambda_1(-1,\operatorname{rec}^{-1}(\sigma))$,
as in \cite[Definition~2.29]{Sato}, and
$\zeta_{\mathfrak f}(s,\sigma)=\zeta_{\mathfrak f}(s,\operatorname{rec}^{-1}(\sigma))$.
We denote by $\widehat{\mathcal L}_2$ the enhanced dilogarithm of
Zagier--Gangl on the rational Bloch group $\mathcal B_2(F)_{\mathbb Q}$, with values in
$\mathbb C/\mathbb Q(2)$.
The definitions and normalizations are given in
Sections~\ref{sec:zeta-def} and~\ref{sec:enhanced}.
The Galois action on Bloch symbols is $\gamma[z]=[\gamma(z)]$.

\begin{theorem}\label{thm:main}
Let $K$ be an imaginary quadratic field and $\mathfrak f$ a nonzero
integral ideal. For every $\sigma\in\operatorname{Gal}(F/K)$ and
every choice of the auxiliary data of
Sections~\ref{sec:eisenstein}--\ref{sec:descent}, formulas
\eqref{eq:bP} and
\eqref{eq:xi} define a class
$\xi_\sigma\in\mathcal B_2(F)_{\mathbb Q}$, which does not depend on
these choices. For all $\sigma,\gamma\in\operatorname{Gal}(F/K)$, we have
$\gamma^{-1}\xi_\sigma=\xi_{\sigma\gamma}$ and
\begin{equation}\label{eq:main}
 \widehat{\mathcal L}_2(\gamma^{-1}\xi_\sigma)=2\pi i\Lambda_1(-1,\sigma\gamma)
       \quad\text{in }\mathbb C/\mathbb Q(2).
\end{equation}
\end{theorem}

With the normalization of $\operatorname{rec}$ above, Theorem~\ref{thm:main}
proves \cite[Conjecture~4.1]{Sato} for $k=1$, with $\xi_{A,2}=\xi_\sigma$
for $A=\operatorname{rec}^{-1}(\sigma)$.

Formula \eqref{eq:xi} defines $\xi_\sigma$ by a
normalized transfer from an auxiliary number field $M$. The restriction of $\xi_\sigma$ to
$M$ is the Galois average of $k_{K,\mathfrak f}b_P$, where $b_P$ is the
explicit finite sum of Bloch symbols \eqref{eq:bP}
and $k_{K,\mathfrak f}$ is the rational constant defined in
\eqref{eq:kKf}.

Taking imaginary parts and using \eqref{eq:ReLambda1} and
Lemma~\ref{lem:inj}, we have the ordinary regulator identity
\[
 \mathcal L_2(\gamma^{-1}\xi_\sigma)=-\pi\zeta_{\mathfrak f}'(-1,\sigma\gamma).
\]
Here $\mathcal L_2$ is the negative of the usual Bloch--Wigner
dilogarithm, as in \cite[\S1.2]{Sato}.
By Borel's theorem, as in the proof of
Lemma~\ref{lem:inj}, these ordinary values for all
$\gamma$ determine the class $\xi_\sigma$. Thus the imaginary part of
\eqref{eq:main} characterizes $\xi_\sigma$, and its real part
is the additional information given by the enhancement.

The classes $\xi_\sigma$ lie in the rational Bloch group, and
\eqref{eq:main} holds modulo $\mathbb Q(2)$. In weight two,
integral refinements involve the extended Bloch group. Building on
Neumann's construction over $\mathbb C$ \cite{Neumann}, Zickert identified
the extended Bloch group of a number field with its indecomposable $K_3$,
including torsion, and gave an explicit formula for the complex regulator
\cite{Zickert}.

Let us briefly explain the idea of the proof. Roughly speaking, we
compute the value of $G_\alpha=(\Delta\partial_\tau+2)F_\alpha$ at the
period ratio $\tau$ of the lattice attached to the ray class in two
different ways. Here $\Delta=\overline\tau-\tau$, and $F_\alpha$ is a cubic
polynomial in $\tau$ plus an explicit $q$-series. On the analytic side,
the Fourier expansion of the partial derivative of the Shintani L-function
at the point $(-2,0)$ is given in \cite[Proposition~3.6]{Sato}, and the
raising identity \cite[Corollary~3.9]{Sato}, which acts by the operator
$\Delta\partial_\tau+2$, carries it to the diagonal point $(-1,-1)$. This
expresses $\Lambda_1(-1,\sigma)$ in terms of $G_\alpha(\tau)$
(Proposition~\ref{prop:fourier}).

On the motivic
side, Doran and Kerr gave explicit higher Chow cycles representing
Beilinson's Eisenstein symbol \cite{Beilinson1986} on the fiber square of
the universal elliptic curve \cite{DK}, and Bloch, Kerr and Vanhove
computed their regulators \cite{BKV}. For a suitable divisor supported on
torsion points, the pairing of the regulator of the corresponding class
with $dz_1\wedge dz_2$ is $F_\alpha$, up to a polynomial in $\tau$ of degree at most two with
coefficients in $\mathbb Q(3)$. When we contract the
Eisenstein class on the square of the CM elliptic curve attached to the
ray class with the difference of the graphs of the complex
multiplications $\bar\rho$ and $\rho$, the same operator
$\Delta\partial_\tau+2$ appears (Remark~\ref{rem:raising}). On both
sides, what we actually compute differs from $G_\alpha$ by a polynomial in
$\tau+\overline\tau$ and $\tau\overline\tau$ with coefficients in
$\mathbb Q(3)$. Since $\tau$ lies in $K$,
these differences are rational multiples of $(2\pi i)^3$, and after the
normalizations in Propositions~\ref{prop:fourier}
and~\ref{prop:regCP} they lie in $\mathbb Q(1)$ and
$\mathbb Q(2)$, respectively. Thus we obtain a higher Chow class of a
number field whose regulator is $2\pi i\Lambda_1(-1,\sigma)$, up to an
explicit rational factor.

The rest of the proof translates this class
into the Bloch group. We show that, under the isomorphism between the
Bloch group and the higher Chow group given by Totaro cycles, the
regulator corresponds to the enhanced dilogarithm in
$\mathbb C/\mathbb Q(2)$. Then we write
the class, after restriction to an auxiliary number field, as an explicit
sum of Bloch symbols, using Bolbachan's reciprocity map and a
quadratic reduction on the elliptic curve adapted from Levin
\cite{Bolbachan,Levin}, and
descend it to the ray class field by complex multiplication and a
normalized transfer.

In Section~\ref{sec:zeta}, we recall the definition of
$\Lambda_1(-1,\sigma)$ and compute its Fourier expansion. In
Section~\ref{sec:eisenstein}, we construct the higher Chow class $C_P$
and compute its regulator. In Section~\ref{sec:bloch}, we
express $C_P$, after restriction to an auxiliary number field, by Bloch
symbols. Finally, in Section~\ref{sec:descent}, we descend to
the ray class field and prove Theorem~\ref{thm:main}.

\section{The enhanced partial zeta value and its analytic expression}
\label{sec:zeta}

In this section, we recall the enhanced partial zeta value
$\Lambda_1(-1,\sigma)$ from \cite{Sato} and, by means of its Fourier
expansion, express it in terms of the function $G_\alpha(\tau)$
(Proposition~\ref{prop:fourier}). Note that the ray
class determines a lattice in $\mathbb C$ together with a torsion point on
the corresponding complex torus.
In Section~\ref{sec:descent}, this pair is realized as a CM
elliptic curve with a torsion point, to which the construction of
Section~\ref{sec:eisenstein} is applied.
Throughout this section, $K$ and $\mathfrak f$ are as in the introduction,
with $K\subset\mathbb C$ under the chosen embedding. Fix $\sigma\in\operatorname{Gal}(F/K)$
and put $A=\operatorname{rec}^{-1}(\sigma)$.

\subsection{Definition, branches, and normalization}
\label{sec:zeta-def}

We first recall the Shintani L-function. Since
$x^{-s_1}\overline x^{-s_2}$ is not single valued in $x$ unless
$s_1-s_2\in\mathbb Z$, it is necessary to work on the universal covering
of $\mathbb C^\times$. Let
$p:\widetilde{\mathbb C^\times}\longrightarrow\mathbb C^\times$
be the universal covering group. We denote by $\widetilde{\log}$ its
logarithm.
For a lift $x$ of a nonzero complex number, the conjugate lift
$\overline x$ is defined by
$\widetilde{\log}(\overline x)=\overline{\widetilde{\log}(x)}$, and
$x^{-s}=\exp(-s\widetilde{\log}(x))$. Fix a lattice
$M\subset\mathbb C$. A cone function is a finite integral
linear combination of characteristic functions of lifted open angular
sectors of angle less than $\pi$ and lifted rays, with all rays
generated by nonzero elements of $M$.
Let $\varepsilon_{\mathrm{deck}}$ be the deck transformation adding
$2\pi i$ to $\widetilde{\log}$, with action on functions
$(\varepsilon_{\mathrm{deck}}f)(x)=f(\varepsilon_{\mathrm{deck}}^{-1}x)$.
A fundamental-domain function is a cone function $D$ satisfying
\[
 \sum_{n\in\mathbb Z}(\varepsilon_{\mathrm{deck}}^nD)(x)=1
 \qquad\text{for every }x\in\widetilde{\mathbb C^\times}.
\]
The sum is finite for each $x$, since a cone function has bounded
angular support. For example, if $\theta$ is an argument of a nonzero
lattice element, then the characteristic function of the half-open domain
$\theta\leq\operatorname{Im}\widetilde{\log}(x)<\theta+2\pi$ is a
fundamental-domain function. When $M$ meets the positive real axis, we
denote by $D(1)$ the one with $\theta=0$.

Let $\phi:M\to\mathbb Z$ be periodic modulo a finite-index sublattice.
For any cone function $C$, define the Shintani L-function by
\begin{equation}
 \mathscr L_M(s_1,s_2;\phi,C)
 =\sum_{x\in p^{-1}(M\setminus\{0\})}
       C(x)\phi(p(x))x^{-s_1}\overline x^{-s_2}.
\end{equation}
The series is absolutely convergent when
$\operatorname{Re}(s_1+s_2)>2$; for its definition and meromorphic
continuation, see \cite[Definition~2.1, Proposition~2.8 and Remark~2.9]{Sato}.
We use the diagonal limiting convention of \cite[\S2.7]{Sato}.
For a fundamental-domain function $D$, set
\begin{equation}
 H_M(\phi,D)
 =\lim_{h\to0}\frac{\partial\mathscr L_M}{\partial s_1}
       (-1+h/2,-1+h/2;\phi,D).
 \label{eq:H}
\end{equation}
In \eqref{eq:H}, we differentiate with
respect to $s_1$ before restricting to the diagonal. The limit exists by
\cite[Lemma~2.12]{Sato}. The vanishing result used below is
\begin{equation}
 \mathscr L_M(-k_1+h/2,-k_2+h/2;\phi,D)=O(h)
 \quad (k_1,k_2\in\mathbb Z_{\geq0},\ k_1+k_2>0),
 \label{eq:vanish}
\end{equation}
which follows from \cite[Propositions~2.8 and~2.11]{Sato}.

Let $\mathfrak a$ be an integral ideal in the class $A^{-1}$ prime to
$\mathfrak f$. Since $\mathfrak a+\mathfrak f=\mathcal O_K$, we can take
$x_0\in\mathfrak a$ with $x_0\equiv1\pmod{\mathfrak f}$. Put
\begin{equation}
 \Lambda=\mathfrak a\mathfrak f,
 \qquad
 w=\#\{u\in\mathcal O_K^\times:u\equiv1\pmod{\mathfrak f}\}.
 \label{eq:Lambda-w}
\end{equation}
Since $\mathfrak a\cap\mathfrak f=\mathfrak a\mathfrak f$, the elements
of $\mathfrak a$ congruent to one modulo $\mathfrak f$ form precisely
$x_0+\Lambda$. Let $\phi_{\mathfrak a}=1_{x_0+\Lambda}$ be its
characteristic function on $\mathfrak a$. Note that $\phi_{\mathfrak a}(0)=1$ exactly when
$\mathfrak f=\mathcal O_K$; the series omits $x=0$ in any case.

The map $x\mapsto(x)\mathfrak a^{-1}$ sends the nonzero elements of
$x_0+\Lambda$ onto the integral ideals in $A$, and each fiber consists of
the $w$ unit multiples in \eqref{eq:Lambda-w}. Each nonzero
$x\in x_0+\Lambda$ satisfies $x\in\mathfrak a$ and
$x\equiv1\pmod{\mathfrak f}$. Then the ideal $(x)\mathfrak a^{-1}$ is
integral, prime to $\mathfrak f$, and in the ray class
$[\mathfrak a]^{-1}=A$. Conversely, if $\mathfrak b$ is an integral
ideal in $A$, then $\mathfrak a\mathfrak b=(x)$ with
$x\equiv1\pmod{\mathfrak f}$, and $\mathfrak b\subset\mathcal O_K$
implies $x\in\mathfrak a$. We denote by
$N(\mathfrak b)=[\mathcal O_K:\mathfrak b]$ the absolute norm of a nonzero
integral ideal. Thus the partial zeta function satisfies, initially for
$\operatorname{Re}s>1$,
\begin{equation}
 \zeta_{\mathfrak f}(s,\sigma)
 =\sum_{\substack{0\ne\mathfrak b\subset\mathcal O_K,\ (\mathfrak b,\mathfrak f)=1\\
                   [\mathfrak b]=A}}N(\mathfrak b)^{-s}
 =\frac{N(\mathfrak a)^s}{w}
       \mathscr L_{\mathfrak a}(s,s;\phi_{\mathfrak a},D).
\end{equation}
The enhanced value is
\begin{equation}
 \Lambda_1(-1,\sigma)
 =\frac{H_{\mathfrak a}(\phi_{\mathfrak a},D)}{wN(\mathfrak a)}
       \pmod{\mathbb Q(1)}.
 \label{eq:Lambda1}
\end{equation}
This is \cite[Definitions~2.25 and~2.28]{Sato} at $k=1$, where these
values are called ray class invariants. Since $\mathbb Q(1)$ is purely
imaginary, the real part of $\Lambda_1(-1,\sigma)$ is well defined, and by
\cite[(1.2) and Theorem~1]{Sato} it satisfies
\begin{equation}\label{eq:ReLambda1}
 \operatorname{Re}\Lambda_1(-1,\sigma)
       =\frac12\zeta_{\mathfrak f}'(-1,\sigma).
\end{equation}

If $D$ and $D'$ are fundamental-domain functions, then
$D-D'=(1-\varepsilon_{\mathrm{deck}})f$ for a cone function $f$
\cite[Proposition~2.4]{Sato}. Hence we have, first in the convergent
half-plane and then by analytic continuation,
\begin{equation}
 \mathscr L_M(s_1,s_2;\phi,D-D')
 =\bigl(1-e^{2\pi i(s_2-s_1)}\bigr)
       \mathscr L_M(s_1,s_2;\phi,f).
\end{equation}
Differentiating with respect to $s_1$ and then restricting to the
diagonal, we have
\begin{equation}
 H_M(\phi,D)-H_M(\phi,D')
 =2\pi i\lim_{h\to0}\mathscr L_M(-1+h/2,-1+h/2;\phi,f).
 \label{eq:branch}
\end{equation}
For $M\subset K$, the last limit is rational by
\cite[Proposition~2.24]{Sato}. The prefactor in
\eqref{eq:Lambda1} is rational, so changing $D$
alters $\Lambda_1(-1,\sigma)$ only by an element of $\mathbb Q(1)$.

Now let $\mathfrak a'=c\mathfrak a$ be another integral ideal in
$A^{-1}$ prime to $\mathfrak f$, with $c\equiv1\pmod{\mathfrak f}$.
Let $\widetilde c$ be a lift of $c$ to $\widetilde{\mathbb C^\times}$.
The weight $\phi_{\mathfrak a'}$ is the transport of $\phi_{\mathfrak a}$
by $c$, since $cx_0\in\mathfrak a'$ and $cx_0\equiv1\pmod{\mathfrak f}$,
and $D_c(x)=D(\widetilde c^{-1}x)$ is a fundamental-domain function for
$\mathfrak a'$.
Then the Shintani L-function is multiplied by
$\widetilde c^{-s_1}\overline{\widetilde c}^{-s_2}$, and suppressing the
auxiliary arguments, we have, by the product rule,
\[
 \frac{\partial\mathscr L_{\mathfrak a'}}{\partial s_1}
 =\widetilde c^{-s_1}\overline{\widetilde c}^{-s_2}
   \left(\frac{\partial\mathscr L_{\mathfrak a}}{\partial s_1}
       -\widetilde{\log}(\widetilde c)\mathscr L_{\mathfrak a}\right).
\]
Along the diagonal tending to $(-1,-1)$, the factor tends to
$|c|^2$, and the last term tends to zero by
\eqref{eq:vanish}. Since
$N(\mathfrak a')=|c|^2N(\mathfrak a)$, it follows that
\begin{align*}
 \frac{H_{\mathfrak a'}(\phi_{\mathfrak a'},D_c)}{N(\mathfrak a')}
 &=\frac{|c|^2H_{\mathfrak a}(\phi_{\mathfrak a},D)}{|c|^2N(\mathfrak a)}\\
 &=\frac{H_{\mathfrak a}(\phi_{\mathfrak a},D)}{N(\mathfrak a)}.
\end{align*}
Together with \eqref{eq:branch}, this shows that the
right-hand side of \eqref{eq:Lambda1} does not
depend on the choices of $D$ and $\mathfrak a$.

\subsection{A Fourier expansion}
\label{sec:fourier}

Now we compute $\Lambda_1(-1,\sigma)$. For this purpose, we first choose a
lattice basis in which the torsion point is a rational multiple of the
first basis vector, so that the Fourier expansion of
\cite[Proposition~3.6]{Sato} takes a simple form.

The point $P=x_0\bmod\Lambda$ of $\mathbb C/\Lambda$ is torsion:
$\mathfrak f x_0\subset\mathfrak a\mathfrak f$. We denote by
$\mathfrak H=\{\tau\in\mathbb C:\operatorname{Im}\tau>0\}$ the upper
half-plane, and take an oriented basis $(\omega_1,\omega_2)$ of $\Lambda$
such that
\begin{equation}
 \tau=\omega_2/\omega_1\in\mathfrak H,
 \qquad P=\alpha\omega_1\bmod\Lambda,
 \qquad 0\leq\alpha<1,\quad\alpha\in\mathbb Q.
\end{equation}
If $P\ne0$, such a basis is obtained as follows: we write a lift of $P$
in an oriented basis with rational coordinates, clear denominators, and
divide the resulting pair of integers by its greatest common divisor; by
B\'ezout's identity, the resulting primitive vector extends to an oriented
integral basis, and reducing its rational coefficient modulo one, we
obtain $\alpha$. If $P=0$, we set $\alpha=0$ and use any oriented basis. Put
\begin{equation}
 a_0=\frac{|\omega_1|^2}{N(\mathfrak a)},
 \qquad \Delta=\overline\tau-\tau,
 \qquad q=e^{2\pi i\tau},\qquad z=e^{2\pi i\alpha}.
 \label{eq:a0}
\end{equation}
Here $a_0\in\mathbb Q$, since $\omega_1\in K$.

Let $N$ be a positive integer such that $r_0=N\alpha$ is integral; we may
enlarge $N$ to match the level used in the algebraic construction. Put
$T=N\tau$ and $b=\omega_1/N$. On
$M_T=\mathbb Z+\mathbb Z T$, define
\begin{equation}
 \phi_T(m+nT)=
 \begin{cases}1,&m\equiv r_0\pmod N,\\0,&m\not\equiv r_0\pmod N.
 \end{cases}
\end{equation}
The support of $\phi_T$ satisfies
\begin{equation}
 b\{m+nT:m\equiv r_0\pmod N\}
 =\alpha\omega_1+\Lambda=x_0+\Lambda.
 \label{eq:support}
\end{equation}
Let $\widetilde b$ be a lift of $b$. We take the cut of argument zero on
$M_T$ and transport it to $bM_T$ by $\widetilde b$. The boundary ray of
the transported cut passes through $\omega_1\in\mathfrak a$, so this cut
is also allowed for $\mathfrak a$. By \eqref{eq:support},
multiplication by $b$ identifies the nonzero terms of the two weighted
series. Then the corresponding Shintani L-functions differ by the
factor $\widetilde b^{-s_1}\overline{\widetilde b}^{-s_2}$, first in the
domain of convergence and then by meromorphic continuation. Repeating
the product-rule calculation above with $b$ in place of $c$, we see that
the term containing $\mathscr L_{M_T}$ tends to zero along the diagonal by
\eqref{eq:vanish}. Thus we have
\begin{equation}
 \Lambda_1(-1,\sigma)
 =\frac{a_0}{wN^2}H_{M_T}(\phi_T,D(1))
       \pmod{\mathbb Q(1)}.
 \label{eq:scaling}
\end{equation}

For $|z|=1$, we use the absolutely convergent series
$\operatorname{Li}_3(z)=\sum_{l\geq1}z^l/l^3$, and define
\begin{align}
 Q_\alpha(\tau)
 &=\sum_{l,j\geq1}\frac{z^l+z^{-l}}{l^3}q^{lj},
 \\
 F_\alpha(\tau)
 &=\frac{(2\pi i)^3\tau^3}{720}
   +\frac{\operatorname{Li}_3(z)+\operatorname{Li}_3(z^{-1})}{2}
   +Q_\alpha(\tau),
 \label{eq:F}\\
 G_\alpha(\tau)
 &=\Delta F_\alpha'(\tau)+2F_\alpha(\tau).
 \label{eq:G}
\end{align}
Primes on $F_\alpha$ and $Q_\alpha$ denote differentiation with
respect to $\tau$, with $\alpha$ fixed. The series for $Q_\alpha$ and
each derivative used here converge locally absolutely: on
$\operatorname{Im}\tau\geq\varepsilon>0$, their summands are bounded
by a polynomial in $l,j$ times $e^{-2\pi\varepsilon lj}$.

\begin{proposition}
\label{prop:fourier}
For every nonzero integral ideal $\mathfrak f$, including
$\mathfrak f=\mathcal O_K$, we have, with the notation above,
\begin{equation}
 \Lambda_1(-1,\sigma)=\frac{a_0}{w(2\pi i)^2}G_\alpha(\tau)
       \quad\text{in }\mathbb C/\mathbb Q(1).
\end{equation}
\end{proposition}

\begin{proof}
We first compute the limiting partial derivative at $(-2,0)$, and then
obtain the limit at $(-1,-1)$ by the raising identity
\eqref{eq:raising}. We apply \cite[Proposition~3.6]{Sato} with $k=2$,
$T\in\mathfrak H$, the lattice $M_T=\mathbb Z+\mathbb Z T$,
the weight $\phi_T$, which is $NM_T$-periodic with $N$ equal to
the integer $N_\phi$ of \cite{Sato}, extended by zero to
$\mathbb Q+\mathbb Q T$, the fundamental-domain function $D(1)$,
and limiting direction $(1/2,1/2)$. Put
\begin{equation}
 J_T=\lim_{h\to0}\frac{\partial\mathscr L_{M_T}}{\partial s_1}
          (-2+h/2,h/2;\phi_T,D(1)).
\end{equation}
The limit exists by \cite[Lemma~2.12]{Sato}, and the formula of
\cite[Proposition~3.6]{Sato} reads
\begin{align}
 \frac{(2\pi i)^2}{2N^2}J_T
 &=\eta_{\phi_T}(3)
 -\frac{\zeta_{\phi_T}(-3)(2\pi i)^3(\overline T-T)^3}{12N^2}
 \nonumber\\
 &\quad+
 \sum_{l\geq1}
 \frac{\displaystyle
   \sum_{m,n=1}^N
    \bigl(\phi_T(m+nT)+\phi_T(-m-nT)\bigr)
       e^{2\pi il(m+nT)/N}}
       {l^3(1-e^{2\pi ilT})}.
 \label{eq:JT}
\end{align}
As $\phi_T$ is supported on $M_T$, the sums over rational coordinates
in that proposition become sums over integers, and the one-dimensional
series there specialize to
\begin{align}
 \zeta_{\phi_T}(s)
 &=\sum_{j\geq1}\frac{1}{Nj^s}
       \sum_{m=1}^N\phi_T(m+jT)
   =\frac{\zeta(s)}{N},
 \\
 \eta_{\phi_T}(3)
 &=\sum_{l\geq1}\frac{1}{l^3}
       \sum_{m=1}^N\phi_T(m)e^{2\pi ilm/N}
   =\operatorname{Li}_3(z).
\end{align}
Here $\zeta$ is the Riemann zeta function, the first equality holds for
$\operatorname{Re}s>1$ and then by meromorphic continuation, and exactly
one residue $m$ contributes to each inner sum. In the numerator of
\eqref{eq:JT}, we extract the common factor
$e^{2\pi ilnT/N}$. Then the two remaining sums over $m$ are $z^l$ and
$z^{-l}$, respectively, and the sum over $n$ satisfies
\begin{equation}
 \frac{\sum_{n=1}^N e^{2\pi ilnT/N}}{1-e^{2\pi ilT}}
 =\frac{q^l}{1-q^l}=\sum_{j\geq1}q^{lj}.
\end{equation}
Then the last term of \eqref{eq:JT} is
$Q_\alpha$. Since
$\zeta(-3)=1/120$ and $\overline T-T=N\Delta$, we obtain
\begin{equation}
 J_T=\frac{2N^2}{(2\pi i)^2}
   \left(\operatorname{Li}_3(z)-\frac{(2\pi i)^3\Delta^3}{1440}
                     +Q_\alpha(\tau)\right).
 \label{eq:JT-value}
\end{equation}

Now we derive the raising identity. For short, we write
$\mathscr L(s_1,s_2)=\mathscr L_{M_T}(s_1,s_2;\phi_T,D(1))$, and denote by
$\partial_T=\partial/\partial T$ the Wirtinger derivative with
$\overline T$ held fixed. For $x=m+nT$, the weight and the sheet selected
by $D(1)$ are locally constant in $T$, since
$\operatorname{Im}x=n\operatorname{Im}T$ and the terms with $n=0$ are
fixed. In the region of absolute convergence, termwise differentiation
is permitted by the bound $|n/x|\leq1/\operatorname{Im}T$, and we have
\[
 (\overline T-T)\partial_T\mathscr L(s_1,s_2)
 =s_1\bigl(\mathscr L(s_1,s_2)-\mathscr L(s_1+1,s_2-1)\bigr).
\]

For $T$ near a fixed point of $\mathfrak H$, we decompose $D(1)$ by the
rays through $1,T,-1,-T$, and for each cone, we choose a fixed rotation placing
its generators strictly in the right half-plane. Then the finite
exponential expressions in the derivation of
\cite[Proposition~2.8]{Sato} (in \S2.4 there) permit fixed,
sufficiently small contour radii: their denominators remain nonzero
on the finite contour pieces, and on the infinite rays the expressions
and their $T$-derivatives are bounded by a polynomial times a common
decreasing exponential. Hence the Taylor coefficients of
$(e^{2\pi i(s_1+s_2)}-1)\mathscr L(s_1,s_2)$ can be differentiated in
$T$, with locally uniform remainders, so that $\partial_T\mathscr L$ is
again meromorphic in $(s_1,s_2)$. By \cite[Proposition~2.11]{Sato}, this expansion starts in degree
two at $(-2,0)$ and $(-1,-1)$. Thus, along the diagonal, the
undifferentiated values are $O(h)$, while the partial derivatives and
their $T$-derivatives have expansions with locally uniform $O(h)$
remainders. In particular, the limit defining $J_T$ commutes with
$\partial_T$.

The differential identity above extends by meromorphic continuation,
and differentiating it with respect to $s_1$, we have
\[
\begin{aligned}
 \bigl((\overline T-T)\partial_T-s_1\bigr)
       \frac{\partial\mathscr L}{\partial s_1}(s_1,s_2)
 &=\mathscr L(s_1,s_2)-\mathscr L(s_1+1,s_2-1)\\
 &\quad-s_1\frac{\partial\mathscr L}{\partial s_1}(s_1+1,s_2-1).
\end{aligned}
\]
At $(s_1,s_2)=(-2+h/2,h/2)$, the two undifferentiated terms tend to zero
as $h\to0$. Taking the limit, we obtain the raising identity of
\cite[Corollary~3.9]{Sato} at $k=1$:
\begin{equation}
 H_{M_T}(\phi_T,D(1))
 =\frac12\left((\overline T-T)\frac{\partial}{\partial T}+2\right)J_T.
 \label{eq:raising}
\end{equation}

Similarly, we put $\partial_\tau=\partial/\partial\tau$, holding
$\overline\tau$ fixed. Since $T=N\tau$ with $N$ and $\alpha$
fixed, we have $(\overline T-T)\partial_T=\Delta\partial_\tau$.
Moreover, we have $(\Delta\partial_\tau+2)\Delta^3=-\Delta^3$, since
$\partial_\tau\Delta=-1$. Using \eqref{eq:scaling},
\eqref{eq:raising}, and \eqref{eq:JT-value}, we obtain
an expression independent of $N$ in $\mathbb C/\mathbb Q(1)$:
\begin{equation}\label{eq:Lambda1-fourier}
\begin{aligned}
 \Lambda_1(-1,\sigma)
 &=\frac{a_0}{2wN^2}
     \bigl((\overline T-T)\partial_T+2\bigr)J_T\\
 &=\frac{a_0}{w(2\pi i)^2}(\Delta\partial_\tau+2)
     \left(\operatorname{Li}_3(z)-\frac{(2\pi i)^3\Delta^3}{1440}
          +Q_\alpha(\tau)\right)\\
 &=\frac{a_0}{w(2\pi i)^2}
     \left(\frac{(2\pi i)^3\Delta^3}{1440}
          +2\operatorname{Li}_3(z)
          +\Delta Q_\alpha'(\tau)+2Q_\alpha(\tau)\right).
\end{aligned}
\end{equation}

Now we compare \eqref{eq:Lambda1-fourier} with
$a_0G_\alpha/(w(2\pi i)^2)$. Put $S=\tau+\overline\tau$ and
$U=\tau\overline\tau$, which are rational since $\tau\in K$. The positive
Fourier terms of the two expressions agree term by term, and after
removing the common factor $a_0/(w(2\pi i)^2)$, the difference of their
$\tau$-independent trilogarithm terms is
$\operatorname{Li}_3(z)-\operatorname{Li}_3(z^{-1})$.
Integrating by parts three times, we see that the Fourier coefficients of
$B_3(x)=x^3-3x^2/2+x/2$ on $[0,1]$ are $-6/(2\pi in)^3$ for $n\ne0$. From
the resulting absolutely convergent Fourier series, we have
\begin{equation}
 \operatorname{Li}_3(e^{2\pi i\alpha})
       -\operatorname{Li}_3(e^{-2\pi i\alpha})
 =-\frac{(2\pi i)^3}{6}B_3(\alpha),\qquad 0\leq\alpha\leq1.
\end{equation}
For the cubic terms, we have
\begin{equation}
 (\Delta\partial_\tau+2)\frac{(2\pi i)^3\tau^3}{720}
 =\frac{(2\pi i)^3}{720}(3\overline\tau\tau^2-\tau^3),
\end{equation}
and we also have the polynomial identity
\begin{equation}
 \Delta^3-2(3\overline\tau\tau^2-\tau^3)=S^3-6US.
\end{equation}
Therefore the difference of the two expressions is
\begin{equation}
 \frac{2\pi ia_0}{w}
       \left(\frac{S^3-6US}{1440}-\frac{B_3(\alpha)}6\right),
\end{equation}
which lies in $\mathbb Q(1)$ since $a_0,S,U,\alpha\in\mathbb Q$.
Thus we complete the proof of the proposition.
\end{proof}

\begin{remark}\label{rem:raising}
The operator $\Delta\partial_\tau+2$ in \eqref{eq:G} comes from
the raising identity \eqref{eq:raising}. In the notation of
\cite[\S3.3]{Sato}, it is $\Psi_{-2}=(\overline\tau-\tau)\partial_{-2}$,
where $\partial_{-2}=\partial_\tau-2/(\tau-\overline\tau)$ is the Maass
raising operator of weight $-2$. The same operator reappears in
Section~\ref{sec:eisenstein} in a completely different way: it comes from
the contraction with the graphs of complex multiplication; see
\eqref{eq:contraction}.
\end{remark}

For the regulator comparison in Section~\ref{sec:descent}, we shall need
the following normalization, which is independent of the ray class.
We denote by $D_K$ the discriminant of $K$, write
$K=\mathbb Q(\sqrt{-d})$ with a squarefree integer $d>0$, and let
$\rho=\sqrt{-d}$ be the square root with positive imaginary part. Since
$\rho\Lambda\subset\Lambda$, there are integers $r,s$ such that
$\rho=r\tau+s$; hence $r>0$ and
$r\operatorname{Im}\tau=\sqrt d$. The covolume of $\mathcal O_K$
is $\sqrt{|D_K|}/2$, and passage to the sublattice
$\Lambda=\mathfrak a\mathfrak f$ multiplies it by its index
$N(\mathfrak a)N(\mathfrak f)$. Computing the same covolume in the
oriented basis $(\omega_1,\omega_1\tau)$, we have
\begin{equation}
 |\omega_1|^2\operatorname{Im}\tau
   =\frac{\sqrt{|D_K|}}2N(\mathfrak a)N(\mathfrak f).
\end{equation}
Dividing by $wN(\mathfrak a)\sqrt d$, we have
\begin{equation}
 \frac{a_0}{wr}
 =\frac{\sqrt{|D_K|}}{2\sqrt d}\frac{N(\mathfrak f)}w
 =:\,k_{K,\mathfrak f}\in\mathbb Q.
 \label{eq:kKf}
\end{equation}
The first factor is $1$ or $1/2$, according to whether
$D_K=-4d$ or $D_K=-d$.

\section{An Eisenstein class and its contraction with the graphs of complex multiplication}
\label{sec:eisenstein}

In this section, we construct a higher Chow class $C_P$ of a number
field whose regulator is $(r/2\pi i)\,G_\alpha(\tau)$
(Proposition~\ref{prop:regCP}). We start from Doran--Kerr's
Eisenstein class on $E\times E$ attached to a torsion divisor; its
regulator was computed by Bloch, Kerr and Vanhove \cite{BKV}, and for a
suitable even divisor $\psi_P$ its pairing with $dz_1\wedge dz_2$ is
$F_\alpha(\tau)$, up to the polynomial ambiguity described in the
introduction; see \eqref{eq:Vpsi}. Contracting this class
with the difference of the graphs of $\bar\rho$ and $\rho$ produces the
operator $\Delta\partial_\tau+2$ of Remark~\ref{rem:raising}, hence
$G_\alpha(\tau)$, and reduces the ambiguity to $\mathbb Q(2)$; see
\eqref{eq:SU}. We carry out the
construction over a number field rather than over $\mathbb C$, since we
need CM reciprocity and a transfer argument in
Section~\ref{sec:descent}.
Throughout this article, higher Chow groups are taken with rational
coefficients.

\subsection{The arithmetic construction}
\label{sec:construction}

Let $E/L$ be an elliptic curve over a number field, let $N>3$, and suppose
that $E[N]\subset E(L)$. We choose a basis of $E[N]$. We identify a
divisor with rational coefficients supported on $E[N]$ with its
coefficient function on $(\mathbb Z/N\mathbb Z)^2$, and define the
convolution of two such divisors by
\begin{equation}
 (\phi*\theta)(c)=\sum_{a+b=c}\phi(a)\theta(b).
\end{equation}
Set
\begin{equation}
 e_N=\frac1{N^2}\sum_{Q\in E[N]}[Q],\qquad
 \delta_N=[0]-e_N.
\end{equation}
Assume that $\deg\phi=0$. Doran--Kerr's construction below takes three
functions as input, while we start from a single divisor $\phi$. To produce such
functions from $\phi$, we use the following factorization. Direct summation shows that
$\phi*e_N=0$ and $e_N*e_N=e_N$. Hence we have
\begin{equation}\label{eq:phi-delta}
 \phi*\delta_N*\delta_N=\phi.
\end{equation}

Since the regulator formula of \cite{BKV} is stated for the universal
family, we carry out the construction over the modular curve of full
level $N$, and then specialize it to $E$. For $N\geq3$, the full level-$N$
moduli problem over $\mathbb Q$
is represented by a smooth affine curve
\cite[Corollary~4.7.2, p.~117]{KatzMazur}; we take $N>3$, the
standing assumption of \cite[\S5.1]{BKV}.
We denote by $\zeta_N$ the primitive $N$th root of unity given by the
Weil pairing of the chosen basis, by $Y_N$ the corresponding component
over $\mathbb Q(\zeta_N)$, and by $\mathcal E\to Y_N$ the universal
elliptic curve. The family is smooth with geometrically integral fibers,
and its total space $\mathcal E$ is a regular surface. All $N$-torsion sections are defined and
pairwise disjoint on $Y_N$.

For distinct torsion sections $A,B$, the divisor $N([A]-[B])$ is
principal on the generic elliptic fiber: under the identification of its
degree-zero Picard group with the elliptic curve, its class is
$N(A-B)=0$. By Hilbert's theorem 90, a defining function can be chosen
over the function field of $Y_N$ itself. We choose a third torsion section
$O'$ and normalize the function $f_{A,B}$ by $f_{A,B}(O')=1$.

On $\mathcal E$, any additional component of $\operatorname{div}(f_{A,B})$
is vertical, and since the fibers are integral, the vertical part is a
sum of whole fibers:
\[
 \operatorname{div}(f_{A,B})
   =N([A]-[B])+\sum_{y\in Y_N^{(1)}}m_y\mathcal E_y,
\]
where $Y_N^{(1)}$ is the set of closed points of $Y_N$.
The section $O'$ is not contained in the support of this Cartier
divisor, so the divisor can be pulled back to $O'$. The horizontal terms
are disjoint from $O'$, while the vertical terms pull back to
$\sum_y m_y[y]$. This pullback is $\operatorname{div}(f_{A,B}(O'))=0$, so
every $m_y$ vanishes. Thus $f_{A,B}$
is a unit on $\mathcal E\setminus\mathcal E[N]$.

Every degree-zero $\mathbb Q$-divisor supported on $\mathcal E[N]$ is a
rational linear combination of the divisors $N([A]-[0])$. Applying this to
each factor of \eqref{eq:phi-delta}, we obtain functions
$f_{ij}$ and coefficients $c_i\in\mathbb Q$ satisfying
\begin{equation}\label{eq:decomposition}
 \phi=\sum_i c_i\,
    \operatorname{div}(f_{i1})*
    \operatorname{div}(f_{i2})*
    \operatorname{div}(f_{i3}).
\end{equation}
We use this decomposition throughout.

Put $\mathcal X=\mathcal E\times_{Y_N}\mathcal E$ and
$\mathcal S=\mathcal E[N]\times_{Y_N}\mathcal E[N]$.
Let $G$ be the group, of order $8N^4$, generated by the translations by
$\mathcal E[N]^2$, the two coordinate inversions, and the interchange of
the two factors, and let $\chi:G\to\{1,-1\}$ be the character that is
trivial on translations and is $-1$ on each inversion and on the
interchange. We use the idempotent
\begin{equation}
 e_\chi=\frac1{8N^4}\sum_{g\in G}\chi(g)g^*.
\end{equation}
For each triple in \eqref{eq:decomposition}, we consider the
ordered graph symbol
\[
 \{f_{i1}(-z_1),\ f_{i2}(z_1-z_2),\ f_{i3}(z_2)\}.
\]
In Doran--Kerr's construction, the closure of its $e_\chi$-projection is
an admissible cubical precycle with zero boundary on
$\mathcal X\setminus\mathcal S$
\cite[\S7.3.4, Proposition~7.17]{DK}. Although \cite{DK} works over
$\mathbb C$, both properties hold for our construction over
$\mathbb Q(\zeta_N)$, by the following direct argument. The divisors of
the $f_{ij}$ are supported on
$\mathcal E[N]$, and if two of $-z_1$, $z_1-z_2$ and $z_2$ lie in
$\mathcal E[N]$, so does the third. Hence, outside $\mathcal S$, at most
one coordinate of the graph symbol has a zero or a pole. Thus the closure
meets all faces properly, and its boundary is supported on the loci
$z_1=A$, $z_2=B$ and
$z_1-z_2=C$ with $A,B,C\in\mathcal E[N]$. These loci are fixed pointwise
by the elements $(z_1,z_2)\mapsto(2A-z_1,z_2)$, $(z_1,2B-z_2)$ and
$(z_2+C,z_1-C)$ of $G$, on which $\chi=-1$. Since
$e_\chi g^*=\chi(g)e_\chi$ for $g\in G$, the projection $e_\chi$ kills
each boundary component. Thus the rational sum with coefficients $c_i$
defines a class on $\mathcal X\setminus\mathcal S$.

The torsion subscheme $\mathcal S$ is smooth of codimension two, and we
have the localization sequence \cite[Properties~17.4(2)]{MVW}
\begin{equation}\label{eq:localization}
 \mathrm{CH}^1(\mathcal S,3)_{\mathbb Q}
 \longrightarrow \mathrm{CH}^3(\mathcal X,3)_{\mathbb Q}
 \longrightarrow \mathrm{CH}^3(\mathcal X\setminus\mathcal S,3)_{\mathbb Q}
 \longrightarrow \mathrm{CH}^1(\mathcal S,2)_{\mathbb Q}.
\end{equation}
We denote by $\mathbb Z(1)_{\mathrm{mot}}$ the motivic Tate complex and
by $\mathbb Q(1)_{\mathrm{mot}}$ its rationalization. From the weight-one
identification $\mathbb Z(1)_{\mathrm{mot}}\simeq\mathbb G_m[-1]$
\cite[Theorem~4.1 and Corollary~4.2]{MVW} and the higher Chow
comparison \cite[Theorem~19.1]{MVW}, we have, for every smooth scheme $W$,
\[
 \mathrm{CH}^1(W,n)_{\mathbb Q}
   =H^{2-n}(W,\mathbb Q(1)_{\mathrm{mot}})
   =H^{1-n}(W,\mathbb G_m\otimes\mathbb Q)=0\qquad(n\ge2).
\]
Hence the outer terms in \eqref{eq:localization} vanish, and the class on
$\mathcal X\setminus\mathcal S$ extends uniquely to
$\mathrm{CH}^3(\mathcal X,3)_{\mathbb Q}$. The extended class remains in the
$\chi$-part, since restriction commutes with $e_\chi$ and is injective.

We pull back at the point of $Y_N$ representing $E/L$ with its level structure. Higher Chow groups
are presheaves with transfers on smooth schemes, in particular
contravariant for morphisms
\cite[Theorem~17.21]{MVW}; concretely, one first moves a representative
so that it meets the fiber and all cubical faces properly. We obtain
\begin{equation}
 Z_\phi\in\mathrm{CH}^3(E^2,3)_{\mathbb Q}^{\chi}.
\end{equation}
Here the superscript $\chi$ denotes the image of $e_\chi$.

\subsection{The Abel--Jacobi regulator}

Next we compute the Abel--Jacobi regulator of the class $Z_\phi$
constructed above, for $N>3$ and an even torsion $\mathbb Q$-divisor
$\phi$ of degree zero, where evenness means $\phi(-Q)=\phi(Q)$. The
regulator of Doran--Kerr's class was computed explicitly in
\cite[Proposition~5.1.1 and (5.1.2)--(5.1.4)]{BKV}, and we have only to
evaluate their formula on our divisor. Note that their hypothesis that
$\widehat\phi$ is even is equivalent to the evenness of $\phi$, by finite
Fourier inversion.

Fix an embedding $L\hookrightarrow\mathbb C$ and write
$E(\mathbb C)=\mathbb C/(\mathbb Z+\tau\mathbb Z)$ with
$\operatorname{Im}\tau>0$. For the calculation, we use the level structure induced by this lattice basis, and take $Y_N$ to be the component
containing this marked curve. We also use the torsion coordinates of
\cite[\S7.2.5, (7.7)]{DK} and the Fourier convention of
\cite[\S5.1]{BKV}:
\begin{equation}
 (x,v)\longleftrightarrow\frac{x\tau+v}{N},\qquad
 \widehat\phi(a,b)
    =\sum_{x,v\bmod N}\phi(x,v)
       \exp\!\left(\frac{2\pi i(av-bx)}N\right).
\end{equation}
Here $\phi(x,v)$ is the coefficient of the point $(x\tau+v)/N$ in
the divisor $\phi$; the arguments of $\widehat\phi$ are read modulo $N$.

Our projector $e_\chi$ is that of \cite[\S7.3.4]{DK}, including the factor
$1/(8N^4)$, and the coefficient function associated with a triple of
functions is the unscaled convolution of their relative divisors
\cite[p.~111, before Theorem~8.6]{DK}. Hence the formula of \cite{BKV}
applies to $Z_\phi$ without further normalization. Moreover, by
\cite[Corollary~9.2]{DK}, the fiberwise regulator depends only on $\phi$,
and not on the functions chosen in \eqref{eq:decomposition}.

As the fiber varies, the universal construction gives an Abel--Jacobi
regulator section $R_\phi$ \cite[\S9.1]{DK}, whose value at $E$ is the
regulator of $Z_\phi$, with target
\[
 H^2(E^2,\mathbb C)/H^2(E^2,\mathbb Q(3)),
\]
since the Hodge filtration term $F^3H^2(E^2,\mathbb C)$ vanishes.
Note that the projector $e_\chi$ retains precisely the three-dimensional
symmetric tensor summand of $H^1(E)\otimes H^1(E)$. Translations act
trivially on Betti cohomology, and the separate inversions kill the two
K\"unneth summands $H^2(E)\otimes H^0(E)$ and $H^0(E)\otimes H^2(E)$. On
$H^1(E)\otimes H^1(E)$, the geometric interchange is minus the tensor
interchange, as both classes have degree one. Thus its minus eigenspace
is $\operatorname{Sym}^2H^1(E)$. In the family, the images
of $e_\chi$ on $H^2(E^2,-)$ form a rank-three local system.

Let $\widetilde R_\phi$ be a holomorphic lift on the upper half-plane
\cite[p.~123]{DK}.
Let $\nabla$ be the Gauss--Manin connection, and
$\langle a,b\rangle=\int_{E^2}a\wedge b$ the Poincar\'e pairing. By
\cite[(5.1.2)]{BKV}, we have
\begin{equation}\label{eq:nablaR}
 \nabla\widetilde R_\phi
   =\kappa_\phi(\tau)[dz_1\wedge dz_2]\,d\tau
\end{equation}
for a scalar function $\kappa_\phi$. Put
$V_\phi=\langle\widetilde R_\phi,[dz_1\wedge dz_2]\rangle$.
In Subsection~\ref{sec:contraction}, we use
\eqref{eq:nablaR} to express the contraction with the graphs of $\rho$ and $\bar\rho$ in terms of $V_\phi$ and its first derivative.
The function $V_\phi$ itself, including its constant term, is given in
\cite[Proposition~5.1.1]{BKV}: modulo $(2\pi i)^3\mathbb Q[\tau]_{\le2}$,
we have, with $q=\exp(2\pi i\tau)$,
\begin{align}\label{eq:Vphi}
 V_\phi(\tau)\equiv{}&
  \frac{2(2\pi i\tau)^3}{(2\pi i)^4}
     \sum_{l\ge1}\frac{\widehat\phi(0,l)}{l^4}
  +\frac1N\sum_{l\ge1}\frac{\sum_{b\bmod N}\widehat\phi(l,b)}{l^3}
   \notag\\
 &+\frac2N\sum_{M\ge1}q^{M/N}
       \sum_{l\mid M}\frac1{l^3}
       \sum_{b\bmod N}
          \exp\!\left(\frac{2\pi ibM}{lN}\right)\widehat\phi(l,b).
\end{align}

\subsection{The divisor attached to \texorpdfstring{$P$}{P}}

Let $P\in E[N]$. We choose an adapted lattice basis as in
\S\ref{sec:fourier} and use the induced level structure, so that $P=\alpha\pmod{\mathbb Z+\tau\mathbb Z}$ with
$\alpha=r_0/N$ and $r_0\in\{0,\ldots,N-1\}$.
Recall $F_\alpha$ and $G_\alpha$ from
\eqref{eq:F} and \eqref{eq:G}. We look for an even divisor $\phi$ of degree zero for which $V_\phi$ is $F_\alpha$ modulo $(2\pi i)^3\mathbb Q[\tau]_{\le2}$. For the
calculation, we extend the right-hand side of
\eqref{eq:Vphi} linearly to all coefficient functions.
The Fourier transform of $([P]+[-P])/2$ is $\cos(2\pi a\alpha)$. Then,
by $\zeta(4)=(2\pi i)^4/1440$, the first term is $(2\pi i)^3\tau^3/720$,
and the second is $\sum_{l\ge1}\cos(2\pi l\alpha)/l^3$. In the last term, the
sum over $b$ is zero unless $N\mid M/l$, and otherwise it is
$N\cos(2\pi l\alpha)$. Substituting $M=lNj$, we see that the last term is $Q_\alpha$. Thus
the linear expression on this divisor is $F_\alpha$.

By the orthogonality of characters, we have $\widehat e_N(a,b)=1$ if
$N\mid a$ and $N\mid b$, and $\widehat e_N(a,b)=0$ otherwise. Then the
cubic term of the expression on $e_N$ is $N^{-4}(2\pi i)^3\tau^3/720$, and
its constant term is $N^{-1}\sum_{N\mid l}l^{-3}=N^{-4}\zeta(3)$. In the
terms with positive powers of $q$, only the indices with $N\mid l$ and
$b=0$ contribute, and reindexing by $l=Nh$ and $M=Nhj$, we have
\begin{align*}
 \frac2N\sum_{M\ge1}q^{M/N}
          \sum_{\substack{l\mid M\\N\mid l}}\frac1{l^3}
 &=\frac2N\sum_{h,j\ge1}\frac{q^{hj}}{(Nh)^3}\\
 &=\frac{2}{N^4}\sum_{h,j\ge1}\frac{q^{hj}}{h^3}\\
 &=N^{-4}Q_0(\tau).
\end{align*}
Here the rearrangement is justified by absolute convergence for $|q|<1$.
Hence the expression on $e_N$ is $N^{-4}F_0$, while its value on $[0]$
is $F_0$.

The divisor $([P]+[-P])/2$ has degree one. Subtracting $e_N$, we obtain a
divisor of degree zero, but this removes $N^{-4}F_0$ from the expression. To restore this term, we add a multiple of the degree-zero
divisor $[0]-e_N$, and put
\begin{equation}
 \psi_P=\frac{[P]+[-P]}2-e_N+\frac{[0]-e_N}{N^4-1}.
\end{equation}
Then the linear expression on $\psi_P$ is
\[
 F_\alpha-N^{-4}F_0+\frac{1-N^{-4}}{N^4-1}F_0=F_\alpha.
\]
Since $\psi_P$ is an even $\mathbb Q$-divisor of degree zero, we have, by
\cite[Proposition~5.1.1]{BKV},
\begin{equation}\label{eq:Vpsi}
 V_{\psi_P}\equiv F_\alpha\pmod{(2\pi i)^3\mathbb Q[\tau]_{\le2}}.
\end{equation}

\subsection{Contraction with the graphs of complex multiplication}
\label{sec:contraction}

Now we use complex multiplication. The class $Z_{\psi_P}$ lives on
$E\times E$. To obtain a class of a number field, we intersect it with a
one-dimensional cycle on $E\times E$ and push the result forward to a
point. As this cycle we take the difference of the graphs of conjugate
CM endomorphisms $\bar\rho$ and $\rho$.

Assume that a nonreal CM endomorphism $\rho$ and its conjugate
$\bar\rho$ are defined over $L$, and identify them with their complex
multipliers under the chosen uniformization. We write
\[
 \rho=r\tau+s,\qquad \rho\tau=u\tau+v,
 \qquad r,s,u,v\in\mathbb Z.
\]
Here $r\ne0$ as $\rho$ is nonreal. From the identity
$r\tau^2+(s-u)\tau-v=0$, we have
\begin{equation}\label{eq:SU}
 S=\tau+\overline\tau=\frac{u-s}{r},\qquad
 U=\tau\overline\tau=-\frac vr,\qquad S,U\in\mathbb Q.
\end{equation}
Let $\Gamma_\rho$ be the graph $z\mapsto(z,\rho z)$ with its complex
orientation. We define the higher Chow class
\begin{equation}\label{eq:CP}
 C_P=-\pi_*\bigl(Z_{\psi_P}\cdot
                    (\Gamma_{\bar\rho}-\Gamma_\rho)\bigr)
             \in\mathrm{CH}^2(L,3)_{\mathbb Q},
 \qquad \pi:E^2\longrightarrow\operatorname{Spec}L.
\end{equation}
Here the product is computed after moving a representative into proper
position; concretely, we pull back to each graph and push forward from
$E$. The construction commutes with finite extension of the base field:
for a finite extension $L'/L$, the class constructed from $E_{L'}$ with
the same auxiliary data is $\operatorname{res}_{L'/L}C_P$, since the fiber
pullback, the graph pullback and the pushforward to the point commute
with the flat base change $\operatorname{Spec}L'\to\operatorname{Spec}L$. We denote by $R:\mathrm{CH}^2(L,3)_{\mathbb Q}\to
\mathbb C/\mathbb Q(2)$ the regulator at the fixed embedding, in the
normalization below.

We write $z=x+\tau y_0$ and put $X=dx$ and $Y=-dy_0$. Their pullbacks to
the two factors give a flat basis of the cohomology summand defined by $e_\chi$:
\[
 X^2=X_1\wedge X_2,\qquad
 XY=\frac{X_1\wedge Y_2+Y_1\wedge X_2}{2},\qquad
 Y^2=Y_1\wedge Y_2.
\]
Using the product complex orientation and $\int_E X\wedge Y=-1$, we
obtain the following matrix for the Poincar\'e pairing in this basis:
\begin{equation}\label{eq:pairing}
 \begin{pmatrix}0&0&-1\\0&1/2&0\\-1&0&0\end{pmatrix}.
\end{equation}
Interchanging the middle two degree-one factors, we have
$\langle X^2,Y^2\rangle=-1$, and squaring the two-form representing
$XY$, we have $\langle XY,XY\rangle=1/2$.
We put $\Omega=(X-\tau Y)^2$ and write a holomorphic lift as
$\widetilde R_{\psi_P}=aX^2+bXY+cY^2$. Then, for
$V=V_{\psi_P}=\langle\widetilde R_{\psi_P},\Omega\rangle$, we have, by
\eqref{eq:pairing},
\[
 V=-a\tau^2-b\tau-c.
\]
As in Section~\ref{sec:zeta}, primes denote derivatives
with respect to $\tau$. By \eqref{eq:nablaR},
there is a scalar function
$\kappa(\tau)$ with
$a'=\kappa$, $b'=-2\tau\kappa$, and $c'=\tau^2\kappa$. Then,
differentiating the expression above, we have
\begin{equation}\label{eq:Vprime}
 V'=-2a\tau-b.
\end{equation}

Multiplication by $\rho$ sends the real coordinates to
$x'=sx+vy_0$, $y_0'=rx+uy_0$, and hence pulls $X$ back to $sX-vY$
and $Y$ to $-rX+uY$. Restricting to its graph, we have
\begin{equation}
 \int_{\Gamma_\rho}X^2=v,\qquad
 \int_{\Gamma_\rho}XY=\frac{s-u}{2},\qquad
 \int_{\Gamma_\rho}Y^2=-r.
\end{equation}
For $\bar\rho$, the coefficients $(r,s,u,v)$ become $(-r,u,s,-v)$, so
the three integrals change sign. Thus we have
\[
 \int_{\Gamma_{\bar\rho}-\Gamma_\rho}\widetilde R_{\psi_P}
   =-2av+b(u-s)+2rc=r(2Ua+Sb+2c).
\]
With $\Delta=\overline\tau-\tau$ as in
\eqref{eq:a0}, substituting
$V=-a\tau^2-b\tau-c$ and \eqref{eq:Vprime}, we have
\begin{align*}
 \Delta V'+2V
 &=\Delta(-2a\tau-b)+2(-a\tau^2-b\tau-c)\\
 &=-2a\tau\overline\tau-b(\tau+\overline\tau)-2c\\
 &=-(2Ua+Sb+2c).
\end{align*}

Thus the graph contraction involves the differential operator appearing in
\eqref{eq:Lambda1-fourier}:
\begin{equation}\label{eq:contraction}
 -\int_{\Gamma_{\bar\rho}-\Gamma_\rho}\widetilde R_{\psi_P}
       =r(\Delta V'+2V).
\end{equation}

In the normalization of \cite[(4.2.2)]{BKV}, a codimension-$p$ cubical
cycle of degree $n$ carries the factor $(2\pi i)^{p-n}$. Graph pullback
preserves $(p,n)=(3,3)$, while pushforward to the point gives $(2,3)$.
Hence the pushforward divides the graph integral by $2\pi i$. Applying
this to \eqref{eq:contraction}, we obtain the regulator of
$C_P$. Here we use that the regulator, Bloch's Abel--Jacobi map computed
by the current of \cite[Theorem~7.1]{KLM}, is compatible with pullback
along the graph embedding and with the proper pushforward to
$\operatorname{Spec}L$, which is integration over $E$.

By \eqref{eq:Vpsi}, we can write $V=F_\alpha+P_2$, where
$P_2=(2\pi i)^3(c_0+c_1\tau+c_2\tau^2)$ with $c_i\in\mathbb Q$. The
contribution of $P_2$ is
\[
 \Delta P_2'+2P_2
       =(2\pi i)^3(2c_0+c_1S+2c_2U).
\]
This belongs to $\mathbb Q(2)$ after multiplication by $r/(2\pi i)$,
by \eqref{eq:SU}. Combining this with \eqref{eq:CP} and
\eqref{eq:contraction}, we obtain the following proposition.

\begin{proposition}\label{prop:regCP}
For every torsion point $P$, including $P=0$, the class $C_P$ defined in
\eqref{eq:CP} satisfies, at the chosen complex embedding,
\begin{equation}
 R(C_P)=\frac{r}{2\pi i}\,G_\alpha(\tau)
                  \quad\text{in }\mathbb C/\mathbb Q(2).
\end{equation}
\end{proposition}

\begin{remark}
When $(E,P)$ is the pointed CM curve attached to a ray class as in
Section~\ref{sec:zeta}, and $\rho=\sqrt{-d}$,
Propositions~\ref{prop:fourier}
and~\ref{prop:regCP}, together with
\eqref{eq:kKf}, show that
$k_{K,\mathfrak f}C_P$ has regulator $2\pi i\Lambda_1(-1,\sigma)$; see
\eqref{eq:regCsigma}.
\end{remark}

\section{The enhanced dilogarithm and the explicit element of the Bloch group}
\label{sec:bloch}

In Section~\ref{sec:eisenstein}, we obtained a higher Chow class with the
expected regulator, while Theorem~\ref{thm:main} is a statement about
Bloch groups. In this section, we translate $C_P$ into the Bloch group.
We first show that, under the isomorphism $\Phi_F$ given by Totaro cycles
(Lemma~\ref{lem:Totaro} and Corollary~\ref{cor:iso}), the
regulator corresponds to Zagier--Gangl's enhanced dilogarithm in
$\mathbb C/\mathbb Q(2)$ (Proposition~\ref{prop:regL2hat}). Then we
use reciprocity on the elliptic curve to write the element of the Bloch group corresponding to the restriction of $C_P$ to an auxiliary number field
as a finite sum of Bloch symbols
(Proposition~\ref{prop:bP}).

\subsection{The enhanced dilogarithm}
\label{sec:enhanced}

We first recall the enhanced dilogarithm of Zagier and Gangl. Note that
its value on a single symbol depends on a choice of
logarithms, and only its sum over an element of the Bloch group is independent of this
choice, modulo $\mathbb Q(2)$. On the Bloch group, its imaginary part is the
Bloch--Wigner dilogarithm, that is, $-\mathcal L_2$ in the notation of the
introduction (Lemma~\ref{lem:inj}).

Let $F\subset\mathbb C$ be a number field containing an imaginary
quadratic field. We put $V_F=F^\times\otimes_{\mathbb Z}\mathbb Q$ and
denote by $\langle z\rangle$ the image of $z\in F^\times$ in $V_F$.
All exterior products in this section are taken over $\mathbb Q$ in
rationalized multiplicative groups.

The rational pre-Bloch group is the $\mathbb Q$-vector space on
symbols $[z]$, $z\in F\setminus\{0,1\}$, modulo the five-term relations
\begin{equation}
 [x]-[y]+[y/x]
 -\left[\frac{1-x^{-1}}{1-y^{-1}}\right]
 +\left[\frac{1-x}{1-y}\right]=0,
 \qquad x,y\ne0,1,\quad x\ne y.
 \label{eq:five-term}
\end{equation}
The map $\delta[z]=\langle z\rangle\wedge\langle1-z\rangle$ vanishes on
these relations, and the rational Bloch group $\mathcal B_2(F)_{\mathbb Q}$
is its kernel on the pre-Bloch group. Thus a
representative $\xi=\sum_i n_i[z_i]$ of an element of the Bloch group
satisfies
\begin{equation}
 \sum_i n_i\langle z_i\rangle\wedge\langle1-z_i\rangle=0.
 \label{eq:delta}
\end{equation}
The group $\mathcal B_2(F)_{\mathbb Q}$ is $B(F)\otimes\mathbb Q$ in the
notation of \cite[\S1]{GZ}.

Let $\ell:V_F\longrightarrow\mathbb C$ be a $\mathbb Q$-linear lift of
the logarithm modulo $\mathbb Q(1)$. Writing $\ell(z)$ for
$\ell(\langle z\rangle)$, we have
\begin{equation}
 \ell(z)-\log z\in \mathbb Q(1).
\end{equation}
Such a lift is obtained by choosing logarithms of the elements of a basis
of $V_F$. Since $\mathbb Q$ is flat over $\mathbb Z$, we have
$V_F\subset V_{\mathbb C}=\mathbb C^\times\otimes_{\mathbb Z}\mathbb Q$, and
extending the basis to $V_{\mathbb C}$, we also obtain a lift on
$\mathbb C^\times$ in the sense of \cite[\S4.1]{SatoThesis}.

For $z\in F\setminus(( -\infty,0]\cup[1,\infty))$,
the enhanced dilogarithm is represented by
\begin{equation}
 \widehat{\mathcal L}_2(z,\ell)
 =\operatorname{Li}_2(z)+\ell(z)\log(1-z)
       -\frac12\ell(z)\ell(1-z).
 \label{eq:L2hat}
\end{equation}
Here $\operatorname{Li}_2$ is continued from its power series without
crossing $[1,\infty)$, and $\log$ denotes the principal logarithm.
This convention follows \cite[\S4.2]{Sato} and Zagier--Gangl; for
the extension of the lift to $\mathbb C^\times$ and the independence of
the values on the Bloch group from the lift, see
\cite[\S4.1, Theorem~4.1]{SatoThesis}.

Following \cite[\S4]{GZ}, we also define the values on the cuts by an
integral. If $\lambda-\log z\in\mathbb Q(1)$, we put
\begin{equation}
 \mathcal F_2(z,\lambda)=\int_0^{\infty}\frac{s-\lambda}{z^{-1}e^s-1}\,ds
       \pmod{\mathbb Q(2)},
\end{equation}
where the contour avoids the poles and agrees with the positive real
axis sufficiently far out. At $s=\log z+2\pi in$, the residue is
$\log z+2\pi in-\lambda\in\mathbb Q(1)$. Hence changing the contour
changes the integral only by an element of $\mathbb Q(2)$. For $|z|<1$, expanding the
denominator and integrating the absolutely convergent series termwise, we
have
\[
 \mathcal F_2(z,\lambda)=\sum_{m\ge1}\frac{z^m}{m^2}
              -\lambda\sum_{m\ge1}\frac{z^m}{m}
 =\operatorname{Li}_2(z)-\lambda\operatorname{Li}_1(z).
\]
The same identity holds modulo $\mathbb Q(2)$ by analytic continuation.
Since $\operatorname{Li}_1(z)=-\log(1-z)$, subtracting
$\ell(z)\ell(1-z)/2$, we obtain \eqref{eq:L2hat}. For arguments
on the cuts, we define the value on a symbol by the same contour:
\[
 \widehat{\mathcal L}_2(z,\ell)=\mathcal F_2(z,\ell(z))-\frac12\ell(z)\ell(1-z)
                      \pmod{\mathbb Q(2)}.
\]
Zagier and Gangl \cite[\S4, Proposition~2]{GZ} show that
its sum on an element of the Bloch group is independent of $\ell$ and of the
representative; they prove the first statement and note that the
second is proved in the same way. For representatives whose arguments
avoid the cuts, the second statement also follows from the proof of
Proposition~\ref{prop:regL2hat} below, which identifies this sum
with the regulator of $\Phi_F(\xi)$.
Thus, for $\xi=\sum_i n_i[z_i]\in \mathcal B_2(F)_{\mathbb Q}$, we define
\[
 \widehat{\mathcal L}_2(\xi)=\sum_i n_i\widehat{\mathcal L}_2(z_i,\ell)
                  \pmod{\mathbb Q(2)}.
\]

\subsection{Comparison with the regulator}

In this subsection, we compare the enhanced dilogarithm with the
regulator. For this purpose, we represent an element of the Bloch group by a sum of Totaro cycles
\cite{Totaro1992,GMS}.
However, this sum is not closed in general, since the condition
\eqref{eq:delta} holds only in $\bigwedge^2V_F$. We correct it by means of bounding chains (Lemma~\ref{lem:Totaro}), and compute the
regulator of the corrected cycle by Wei\ss schuh's path regulator.

\subsubsection{Totaro cycles}

Put $\square=\mathbb P^1\setminus\{1\}$. The rational cubical cycle
group in degree $n$ consists of codimension-$p$ cycles on
$\operatorname{Spec}F\times\square^n$ meeting all faces properly,
modulo cycles pulled back along coordinate projections. A face is obtained
by setting any collection of coordinates to zero or infinity.
Let $\operatorname{Alt}_n$ be the signed average over the $2^n n!$
transformations obtained by permuting the coordinates and inverting
any subset of them, where the sign is the sign of the permutation, with
a factor of $-1$ for each inverted coordinate
\cite[\S II, p.~2]{GMS}. We denote by $C^p(F,n)$ the image of
$\operatorname{Alt}_n$, with the differential
\[
 d=\sum_{j=1}^n(-1)^{j-1}
       (\partial_{j,0}-\partial_{j,\infty}).
\]
Here $\partial_{j,\varepsilon}$ is the intersection with the face whose
$j$th coordinate is $\varepsilon\in\{0,\infty\}$. The product is the
external product followed by the normalized alternating projector, and it
is graded commutative for the cohomological degree $2p-n$. The complex
$C^p(F,\bullet)$ computes $\mathrm{CH}^p(F,n)_{\mathbb Q}$
\cite[Theorems~4.7 and~4.11]{Levine1994}.

Let $s(z)\in C^1(F,1)$ be the alternating point chain for
$z\in F^\times\setminus\{1\}$, and put $s(1)=0$.
From the weight-one calculation
\cite[Theorem~4.1, Corollary~4.2, and Theorem~19.1]{MVW}, we have
\begin{equation}
 C^1(F,1)/dC^1(F,2)\simeq V_F,
 \qquad [s(z)]\longmapsto\langle z\rangle,
 \qquad \mathrm{CH}^1(F,2)_{\mathbb Q}=0.
 \label{eq:weight-one}
\end{equation}
To lift \eqref{eq:delta} to the level of cycles, we choose point
chains representing the elements of $V_F$ linearly, and compare them with $s(z)$ by means of bounding chains.

The $F$-rational point chains already surject onto $V_F$, and we choose a
linear section $u:V_F\longrightarrow C^1(F,1)$ of the quotient map
$C^1(F,1)\longrightarrow V_F$ with image in their span.
For each argument needed below, we choose $h(z)\in C^1(F,2)$ such that
\begin{equation}
 dh(z)=s(z)-u(\langle z\rangle).
 \label{eq:filling}
\end{equation}
We write $u(z)$ for $u(\langle z\rangle)$.

For $z\ne0,1$, brackets in the following formula denote
the ordered parametrized curve before alternation:
\[
 T(z)=\operatorname{Alt}_3[x,1-z/x,1-x]\in C^2(F,3).
\]
Before alternation, the curve has its third, first, and second
coordinates equal to one at $x=0,1,\infty$, respectively, so these
points lie outside the cube. The remaining face is the second zero face
$x=z$, with remaining coordinates $z$ and $1-z$, and its sign is minus.
Thus we have
\begin{equation}
 dT(z)=-s(z)s(1-z).
 \label{eq:dT}
\end{equation}

\begin{lemma}\label{lem:Totaro}
If $\xi=\sum_i n_i[z_i]$ satisfies \eqref{eq:delta}, then
\begin{equation}
 C_\xi=\sum_i n_i\bigl(
 T(z_i)+h(z_i)s(1-z_i)-u(z_i)h(1-z_i)\bigr)
 \label{eq:Cxi}
\end{equation}
is closed, its higher Chow class depends only on $\xi$, and
$\xi\mapsto[C_\xi]$ defines a homomorphism
\[
 \Phi_F:\mathcal B_2(F)_{\mathbb Q}\longrightarrow\mathrm{CH}^2(F,3)_{\mathbb Q},
\]
which is the negative of the Gangl--M\"uller-Stach Totaro-symbol map.
\end{lemma}

\begin{proof}
Since the cohomological degrees of $s$ and $u$ are one and that of $h$ is
zero, we have, by the Leibniz rule together with \eqref{eq:filling} and
\eqref{eq:dT},
\begin{align*}
 d\bigl(T(z)+h(z)s(1-z)-u(z)h(1-z)\bigr)
 &=-s(z)s(1-z)+(s(z)-u(z))s(1-z)\\
 &\quad+u(z)(s(1-z)-u(1-z))\\
 &=-u(z)u(1-z).
\end{align*}
By the linearity of $u$ and the alternating property of the
multiplication, the sum of these boundaries is the image of $-\delta\xi$,
which is zero.

To identify the class of $C_\xi$, we use the acyclic subcomplex of
\cite[\S II.I, pp.~2--3]{GMS}. In cubical degree three, it contains the
products $C^1(F,1)\,C^1(F,2)$, and in cubical degree two, it contains
$C^1(F,1)\,dC^1(F,2)$. The quotient map is a quasi-isomorphism, and both
correction terms in \eqref{eq:Cxi} lie in $C^1(F,1)\,C^1(F,2)$ and hence
vanish in the quotient.
Gangl--M\"uller-Stach's map on kernels factors through
$\mathcal B_2(F)_{\mathbb Q}$
\cite[\S II.I, paragraph following Corollary~2, p.~7]{GMS} and is
represented by $C_z=\operatorname{Alt}_3[x,1-x,1-z/x]$, and swapping its
last two coordinates, we have $T(z)=-C_z$. Hence $[C_\xi]$ is the
negative of the image of $\xi$ under Gangl--M\"uller-Stach's map. This
proves the stated sign, and shows that the class of $C_\xi$ does not
depend on the choice of the section, the bounding chains, and the representative of $\xi$. The
quotient and its symbol formula commute with field isomorphisms, and so
does $\Phi_F$. Thus we complete the proof of the lemma.
\end{proof}

\subsubsection{Representatives in good real position}

We may assume that the arguments of a Bloch representative avoid
$(-\infty,0]$ and $[1,\infty)$. If a real argument $a\ne0,1$
occurs, we choose a nonreal $y\in F$, which exists as $F$ contains an
imaginary quadratic field, and solve \eqref{eq:five-term} with $x=a$ for
$[a]$. Then the other arguments
\[
 y,\qquad y/a,\qquad
 \frac{(a-1)y}{a(y-1)},\qquad \frac{1-a}{1-y}
\]
are the values at $y$ of nonconstant real M\"obius transformations, so
they are nonreal. Applying this replacement to the finitely many real
arguments, we obtain the required representative, and the enhanced value
does not change under the five-term relations.

We may also choose the section $u$ so that its point chains avoid the
negative real axis. For this, we select a basis of $V_F$ from the spanning
set $\{\langle b\rangle:b\in F^\times\}$, and if a selected representative
$b$ is negative real, we use $s(b^2)/2$ instead of $s(b)$.

We say that a chain is in good real position if it meets properly all
intersections of the coordinate logarithm cuts, on every cubical face.
For the ordered Totaro curve, the three open negative-real cuts pull back
to the real arcs $x<0$, $x=tz$ with $0<t<1$, and $x>1$, respectively.
The first and third are disjoint. A common point of the first two would
force $z<0$, and one of the last two would force $z>1$; both are excluded
by the choice of the representative. Hence all pairwise and triple
intersections of the cuts are empty. The endpoints $x=0,1,\infty$ lie outside the cube,
and the only remaining cubical face is $x=z$, where the other coordinates
$z$ and $1-z$ avoid the cuts. Thus the ordered curve is in good real
position.

After extension to $\mathbb C$, the bounding chains can be put in good real position while keeping their boundary fixed. The cycle
$b_z=s(z)-u(z)$ is already in good real position, and its class is zero in
the ordinary rational cubical complex by \eqref{eq:filling}. The inclusion
of the subcomplex of chains in good real position is a rational
quasi-isomorphism by the real moving lemma of Kerr and Lewis
\cite[Lemma~8.14]{KerrLewis}; for the conditions on all intersections of
the cuts used here, see \cite[\S2, Proposition~1]{Weissschuh}. Injectivity on homology then yields a bounding chain $h_0$ in good real position with
$dh_0=b_z$.

We apply the full alternating projector defined above to $h_0$.
Coordinate permutations and inversions preserve the cubical faces and the
negative-real cuts, and thus all the conditions of good real position.
Put $h=\operatorname{Alt}_2h_0$. Since the projectors commute with the
differential and $b_z$ is already fully alternating, we have
\[
 dh=\operatorname{Alt}_1b_z=b_z.
\]
Thus $h$ is a bounding chain in good real position in the required complex. The difference between $h$ and the original bounding chain is a closed weight-one
chain in degree two, which is a boundary by \eqref{eq:weight-one}, now
over $\mathbb C$. By the Leibniz rule, multiplication by the closed point
chains in \eqref{eq:Cxi} preserves boundaries. Thus the
replacement does not change the image of $C_\xi$ in
$\mathrm{CH}^2(\mathbb C,3)_{\mathbb Q}$. The product corrections in
\eqref{eq:Cxi} are also in good real position: their constant
coordinates avoid the cuts and the cubical faces, and the remaining
intersection conditions reduce to those for the bounding chains.

\subsubsection{Calculation of the path regulator}
Let $v$ be the auxiliary path parameter in Wei\ss schuh's commutative
regulator $r_P$, which is written $x$ in \cite{Weissschuh}; we use
rational coefficients throughout. At a point, the target of $r_P$ is
\[
 P_{\mathcal D}(\operatorname{Spec}\mathbb C,\mathbb Q(p))
 =\{w\in\mathbb C[v,dv]:w(0)\in\mathbb Q(p),\ w(1)=0\},
 \qquad p>0.
\]
Here $\mathbb C[v,dv]$ is the algebra of polynomial differential forms. For
degree-one paths, the endpoint condition is automatic, and the
evaluation in Deligne cohomology is the integration from zero to one
modulo $\mathbb Q(p)$; see \cite[\S\S5.1 and 6.1]{Weissschuh}.
The regulator is multiplicative on the alternating complex.

Coordinate permutations contribute their ordinary signs, and coordinate
inversion changes the signs of $\log z$ and $d\log z$.
The negative-real cut current also changes sign, since inversion
interchanges the endpoints zero and infinity. Thus we have
$r_P(\operatorname{Alt}_n Z)=r_P(Z)$ for the chains $Z$ in good real
position used above.

The section $u$ induces a logarithm lift $\ell_u$, characterized by
$r_P(u(z))=\ell_u(z)\,dv$. For one argument $z$, we put
\[
 a=\log z,\qquad b=\log(1-z),\qquad
 a_u=\ell_u(z),\qquad b_u=\ell_u(1-z).
\]
The weight-one regulator gives
\[
 r_P(s(z))=a\,dv,\qquad r_P(u(z))=a_u\,dv.
\]
From \eqref{eq:filling}, we have
$d r_P(h(z))=(a-a_u)dv$. Since a degree-zero path is a polynomial
in $v$ and $r_P(h(z))(1)=0$, we obtain
\begin{equation}
 r_P(h(z))=(a-a_u)(v-1).
 \label{eq:rPh}
\end{equation}
The value of this path at zero is $a_u-a\in\mathbb Q(1)$.

For the coordinate order of $T(z)$, we have, by
\cite[\S6.3, pp.~20--21]{Weissschuh},
\[
 r_P(T(z))=v(1-v)\bigl(6\operatorname{Li}_2(z)+3ab\bigr)\,dv.
\]
The integral of this form is $\operatorname{Li}_2(z)+ab/2$. By
\eqref{eq:rPh}, the correction terms contribute
\begin{align*}
 \int_0^1 r_P(h(z)s(1-z))
 &=(a-a_u)b\int_0^1(v-1)\,dv=-\frac12(a-a_u)b,\\
 \int_0^1 r_P(-u(z)h(1-z))
 &=-a_u(b-b_u)\int_0^1(v-1)\,dv=\frac12a_u(b-b_u).
\end{align*}
Adding these three expressions, we have
\begin{equation}
 \operatorname{Li}_2(z)+\frac12ab
   -\frac12(a-a_u)b+\frac12a_u(b-b_u)
 =\operatorname{Li}_2(z)+a_ub-\frac12a_ub_u.
 \label{eq:rPsum}
\end{equation}

\begin{proposition}\label{prop:regL2hat}
For every $\xi\in \mathcal B_2(F)_{\mathbb Q}$, the
Kerr--Lewis--M\"uller-Stach (KLM) regulator
\cite[\S\S5.5 and 5.7]{KLM}, in the normalization of Wei\ss{}schuh
used here, satisfies
\begin{equation}
 R\bigl(\Phi_F(\xi)\bigr)=\widehat{\mathcal L}_2(\xi)
       \quad\text{in }\mathbb C/\mathbb Q(2).
\end{equation}
\end{proposition}

\begin{proof}
Summing \eqref{eq:rPsum} over the corrected closed
representative, we see that the path regulator of $C_\xi$ is the sum of
\eqref{eq:L2hat} for the lift $\ell_u$ induced by $u$. We show
that the sum does not depend on the lift.
Let $\ell'=\ell_u+2\pi i\mu$ be another lift, where $\mu:V_F\to\mathbb Q$ is
linear. Then, for one symbol, the difference of the two values of
\eqref{eq:L2hat} is exactly
\[
 2\pi i\mu(z)b-\pi i\bigl(\mu(z)b_u+a_u\mu(1-z)\bigr)
       -\frac{(2\pi i)^2}{2}\mu(z)\mu(1-z).
\]
Since $b-b_u\in\mathbb Q(1)$, we have, for the Bloch representative
$\xi=\sum_i n_i[z_i]$, in $\mathbb C/\mathbb Q(2)$,
\begin{align*}
 \sum_i n_i\bigl(\widehat{\mathcal L}_2(z_i,\ell')
                    -\widehat{\mathcal L}_2(z_i,\ell_u)\bigr)
 &=\pi i\sum_i n_i\bigl(\mu(z_i)\ell_u(1-z_i)
                         -\ell_u(z_i)\mu(1-z_i)\bigr)\\
 &=\pi i(\mu\wedge\ell_u)(\delta\xi)\\
 &=0.
\end{align*}
Here $(\mu\wedge\ell_u)(x\wedge y)=\mu(x)\ell_u(y)-\ell_u(x)\mu(y)$,
and the last equality follows from \eqref{eq:delta}. Thus the
sum does not depend on the choice of the logarithm lift.

It remains to compare this path regulator with the regulator $R$ of
Section~\ref{sec:eisenstein}. In cube coordinates $z_1,z_2,z_3$, moving the
first cut current past
$d\log z_3$ changes the sign of the middle term in the three-coordinate
current of \cite[\S1.2, p.~26]{DK}, which is also the current of
$\{f_1,f_2,f_3\}$ displayed in \cite[\S4.2]{BKV}. The ordered intersection of the two cut
currents is the negative of the current of integration over the
Cartesian product of the two oriented cuts. With these two sign changes,
we obtain the formula in \cite[p.~12]{Weissschuh}.
The path and KLM current regulators agree on the corrected closed cycle by
\cite[\S5.1, Lemmas~4--6; \S7, Corollary~8]{Weissschuh}.
Their normalization at a point is $(2\pi i)^{2-3}=(2\pi i)^{-1}$
\cite[\S\S6.1 and 6.3]{Weissschuh}, as in
Subsection~\ref{sec:contraction}.
\end{proof}

\subsubsection{Imaginary parts, injectivity, and the isomorphism \texorpdfstring{$\Phi_F$}{Phi\_F}}
Following \cite[\S1.2]{Sato}, we define
$\mathcal L_2:\mathbb C\setminus\{0,1\}\to\mathbb R$ by
\[
 \mathcal L_2(z)=-\operatorname{Im}\operatorname{Li}_2(z)
                  -\log|z|\arg(1-z)
\]
away from the branch cuts, using principal branches, and extend it
continuously across the cuts. Thus $\mathcal L_2$ is the negative of the
usual Bloch--Wigner dilogarithm. For an element of the Bloch group represented by
$\xi=\sum_i n_i[z_i]$, we put $\mathcal L_2(\xi)=\sum_i n_i\mathcal L_2(z_i)$.

\begin{lemma}\label{lem:inj}
For every $\xi\in\mathcal B_2(F)_{\mathbb Q}$,
\begin{equation}\label{eq:ImL2hat}
 \operatorname{Im}\widehat{\mathcal L}_2(\xi)=-\mathcal L_2(\xi).
\end{equation}
Let $K\subset F$ be an imaginary quadratic field with a fixed complex
embedding. The enhanced values of an element of $\mathcal B_2(F)_{\mathbb Q}$ at all embeddings $F\hookrightarrow\mathbb C$ extending this embedding of $K$
determine the class uniquely.
\end{lemma}

\begin{proof}
Let $\ell=r+i\vartheta$ with $r(z)=\log|z|$. Expanding
\eqref{eq:L2hat}, we have
\[
 \operatorname{Im}\widehat{\mathcal L}_2(z,\ell)
  =-\mathcal L_2(z)+\frac12\bigl(\vartheta(z)r(1-z)-r(z)\vartheta(1-z)\bigr).
\]
Since $r$ and $\vartheta$ are $\mathbb Q$-linear on $V_F$, the additional terms sum to zero on an element of the Bloch group by \eqref{eq:delta},
applied to their alternating pairing. Using representatives off the
cuts, which is possible by the five-term argument above, we obtain
$\operatorname{Im}\widehat{\mathcal L}_2(\xi)=-\mathcal L_2(\xi)$. Note
that this imaginary part is well defined, since $\mathbb Q(2)$ is real.

Since $K$ is imaginary quadratic, each conjugate pair of embeddings
$F\hookrightarrow\mathbb C$ has exactly one member extending the fixed
embedding of $K$. Hence the ordinary values at these embeddings form
the full ordinary regulator vector. This regulator is injective on
$\mathcal B_2(F)_{\mathbb Q}$ by Borel's theorem \cite{Borel1977}
and its Bloch-group interpretation \cite[Theorem~5.2]{Suslin1991}; see also
\cite[\S1]{GZ}. Therefore, if two classes have the same enhanced values,
then they have the same ordinary values by
\eqref{eq:ImL2hat}, and their difference is zero.
\end{proof}

Gangl and M\"uller-Stach note that their map realizes the well-known
isomorphism $\mathcal B_2(F)_{\mathbb Q}\simeq\mathrm{CH}^2(F,3)_{\mathbb Q}$
\cite[\S II.I, p.~7]{GMS}. Since we use the bijectivity of $\Phi$ below,
we include a short proof for number fields.

\begin{corollary}\label{cor:iso}
Let $F\subset\mathbb C$ be a number field containing an imaginary
quadratic field. Then $\Phi_F$ is an isomorphism. Moreover, $\Phi_k$
is an isomorphism for an algebraic closure $k$ of $\mathbb Q$.
\end{corollary}

\begin{proof}
Let $\xi\in\mathcal B_2(F)_{\mathbb Q}$ satisfy $\Phi_F(\xi)=0$. For every
embedding $\iota:F\hookrightarrow\mathbb C$, we have
$\Phi_{\iota(F)}(\iota\xi)=\iota_*\Phi_F(\xi)=0$, since $\Phi$ commutes with
field isomorphisms. Hence all enhanced values of $\xi$ vanish by
Proposition~\ref{prop:regL2hat}, and $\xi=0$ by
Lemma~\ref{lem:inj}. Thus $\Phi_F$ is injective. On the
other hand, $\mathrm{CH}^2(F,3)_{\mathbb Q}$ is isomorphic to the weight-two
Adams eigenspace of $K_3(F)\otimes\mathbb Q$ \cite{Bloch1986,Levine1994},
whose dimension is at most the number $r_2(F)$ of complex places of $F$
by Borel's theorem \cite{Borel1974}.
Since $\dim_{\mathbb Q}\mathcal B_2(F)_{\mathbb Q}=r_2(F)$ by the theorems of
Borel and Suslin \cite[Theorem~5.2]{Suslin1991} (see also \cite[\S1]{GZ}),
$\Phi_F$ is bijective.

For the last assertion, note that every element of
$\mathcal B_2(k)_{\mathbb Q}$ or $\mathrm{CH}^2(k,3)_{\mathbb Q}$ is the
restriction of an element over a number field containing a fixed
imaginary quadratic field, and that $\Phi$ commutes with restriction.
Moreover, restriction to $k$ is injective on rational higher Chow groups,
since a relation over $k$ is defined over a finite extension $F'$ of $F$,
and corestriction after restriction from $F$ to $F'$ is multiplication by
$[F':F]$ \cite[Exercise~1.11 and Theorem~17.21]{MVW}. Hence the
bijectivity of $\Phi_F$ for all such $F$ implies that of $\Phi_k$.
\end{proof}

\subsection{Reciprocity on curves and the evaluation of \texorpdfstring{$h_E$}{h\_E}}

Next we write the element of the Bloch group corresponding to $C_P$ explicitly. For this
purpose, we use Bolbachan's reciprocity map, which induces, up to sign, the
inverse of the isomorphism $\Phi$ of Corollary~\ref{cor:iso}
(Lemma~\ref{lem:Psi}). On an elliptic curve, this map
can be evaluated by a quadratic reduction, which reduces the
computation to $\mathbb P^1$
(Lemma~\ref{lem:RE}).

Fix an algebraic closure $k=\overline L$ of a number field $L$.
For a smooth projective integral curve $C/k$ and nonzero rational functions
$u,v$ on $C$, we define
\begin{equation}\label{eq:SC}
 S_C(u,v)=\sum_{Q\in C}\operatorname{ord}_Q(v)[u(Q)],
 \qquad [0]=[1]=[\infty]=0.
\end{equation}
Here the sum is over the geometric zeros and poles of $v$, and we write
$S_E$ when $C=E$. For the same curve $C/k$, let $h_C=-H_{k(C)}$ be
Bolbachan's reciprocity map, which is defined on
$\bigwedge^3(k(C)^\times\otimes_{\mathbb Z}\mathbb Q)$
and takes values in the rational pre-Bloch group of $k$; as in
\cite{Bolbachan}, the base field $k$ is algebraically closed of
characteristic zero.
From \cite[Definition~1.6, the definition of $\operatorname{RecMaps}$ in
\S1.2, and Theorem~1.8]{Bolbachan}, we have
\begin{align}\label{eq:hC}
 h_C(c\wedge f\wedge g)&=0 \quad(c\in k^\times),\notag\\
 h_C(u\wedge(1-u)\wedge v)&=S_C(u,v)\quad(u\ne1),\notag\\
 h_D(\pi^*w)&=(\deg\pi)h_C(w)
       \quad(\pi:D\longrightarrow C\text{ finite}).
\end{align}
Here the second identity follows from the commutative diagram in
\cite[Definition~1.6(i)]{Bolbachan} and the residue formula in
\cite[Proposition~1.4]{Bolbachan}. With this sign convention, $\delta h_C$ is the sum of the tame symbols at all points of $C$. For an integral admissible curve
in the three-cube, we apply the reciprocity map to its ordered coordinate
functions on the smooth projective normalization, and denote by $\Psi$
the resulting linear operation. We also use the Totaro comparison
$\Phi_k$ over $k$, defined as in Lemma~\ref{lem:Totaro}; it is
an isomorphism by Corollary~\ref{cor:iso}.

\begin{lemma}\label{lem:Psi}
The operation $\Psi$ induces a map from
$\mathrm{CH}^2(k,3)_{\mathbb Q}$ to $\mathcal B_2(k)_{\mathbb Q}$, and
\begin{equation}\label{eq:PhiPsi}
 \Phi_k^{-1}=-\Psi
 \quad\text{on }\mathrm{CH}^2(k,3)_{\mathbb Q}.
\end{equation}
\end{lemma}

\begin{proof}
We first show that $\Psi$ takes values in the Bloch group $\mathcal B_2(k)_{\mathbb Q}$. Note that
$\Psi$ respects cubical alternation, and that degenerate and
constant-coordinate cycles have value zero. The degree formula in
\eqref{eq:hC} shows that a finite parametrization gives
exactly the multiplicity prescribed by proper pushforward. On the projective
normalization, a point outside the cube contributes zero, since one
coordinate is a unit with residue one. Inside the cube, admissibility
permits at most one coordinate to have a zero or pole, and if all
coordinates are units, the tame residue is zero. The remaining residues
are the ordered face wedges, with the signs of the cubical differential
and the face multiplicities obtained by summing the orders over the
branches, and these contributions cancel for a closed cycle. Since $\delta h_C$ is the sum of the tame symbols, we see that $\Psi$ takes values in the Bloch group.

Now we show that $\Psi$ kills boundaries. Let $W$ be an admissible
surface in the four-cube. We take the proper closure of $W$ in
$(\mathbb P^1)^4$, normalize it and resolve it; the coordinate functions
$t_1,\ldots,t_4$ remain morphisms to $\mathbb P^1$. Every divisorial
valuation of the function field of $W$ has a center on this proper
surface. We classify the valuations by the center on the closure of $W$,
counting an exceptional curve of the resolution with the point of the
closure that it lies over; the coordinates that are units along such a
curve are constant, with their values at that point. We apply the surface
reciprocity law \cite[Corollary~1.13]{Bolbachan} to their exterior
product. A divisorial valuation with center outside the cube has a
coordinate with residue $1$. Hence the tame residue of the exterior
product is zero. For a divisorial valuation centered at a point inside
the cube, admissibility excludes intersections with three coordinate
faces. Then at least two coordinates are units with nonzero constant
residues. Every term of the tame residue retains these two factors, so
its image under the reciprocity map is zero. For a divisorial valuation
centered at a curve inside the cube, at most one coordinate has a zero or
pole at the generic point of that curve, and if all four coordinates are
units, the residue is zero. If $t_j$ has the zero or pole along a branch
valuation $\nu$, then we have, by the tame-symbol convention of
\cite[Proposition~1.4]{Bolbachan},
\[
 \partial_\nu(t_1\wedge t_2\wedge t_3\wedge t_4)
 =(-1)^{j-1}\nu(t_j)\,
   \overline{t_1}\wedge\cdots\wedge\widehat{\overline{t_j}}
                   \wedge\cdots\wedge\overline{t_4},
\]
where the bars denote residues of the unit coordinates. The positive
order at zero and the negative order at infinity give the sign in
$(-1)^{j-1}(\partial_{j,0}-\partial_{j,\infty})W$. The degree formula
for $h$ weights the branches by their residue-field degrees, so their sum
gives the face multiplicity in $dW$. Thus, by the surface reciprocity
law, we have $\Psi(dW)=0$.

Finally, we determine the sign. For the ordered Totaro curve used in
$\Phi_k$, we have
\[
 x\wedge(1-a/x)\wedge(1-x)
       =-x\wedge(1-x)\wedge((x-a)/x).
\]
Thus $\Psi(T(a))=-[a]$ by the second identity in \eqref{eq:hC}. Since the product
corrections in the definition of $\Phi_k$ have a constant coordinate, they
contribute zero. Hence $\Psi\Phi_k=-1$ on the Bloch group, and
\eqref{eq:PhiPsi} follows from Corollary~\ref{cor:iso}.
Thus we complete the proof of the lemma.
\end{proof}

Let $E/k$ have a short Weierstrass equation
$y^2=x^3+a_E x+b_E$. To evaluate $h_E$, we use a quadratic reduction
based on an identity of \cite[\S6.2, Lemma~6.2.1]{Levin}. We separate off
the Steinberg terms, average the
remaining exterior products under an involution, and evaluate the
resulting invariant products on the quotient $\mathbb P^1$ using the
degree formula.

\subsubsection{Choice of involution}
Let $f_1,f_2,f_3$ be nonzero rational functions on $E$. We define an
element $\mathcal R_E(f_1,f_2,f_3)$ of the rational pre-Bloch group as
follows. If one of $f_1,f_2,f_3,f_1/f_2,f_1/f_3,f_2/f_3$ is constant, we
put $\mathcal R_E(f_1,f_2,f_3)=0$. Otherwise, let $U$ be an algebraic
point such that none of these six functions is invariant under the
involution $\iota(Q)=2U-Q$. Such $U$ exists. In fact, the translations
preserving a nonconstant rational function form a finite group: they act
faithfully on any finite fiber outside its branch locus. If two
involutions $Q\mapsto2U-Q$ and $Q\mapsto2U'-Q$ preserve the function, then
their composition is the translation by $2(U-U')$. Since this
translation lies in that finite group and $[2]:E\to E$ has finite fibers,
only finitely many $U$ are excluded for each function. In
particular, a single $U$ can be chosen for any finite
list of triples.

For a function $f$ on $E$, we write $f^\iota=f\circ\iota$. Put
\begin{equation}
\begin{array}{lll}
 A=(f_2^{\iota}-f_2)/(f_1f_2^{\iota}-f_2f_1^{\iota}),& B=(f_1-f_1^{\iota})/(f_1f_2^{\iota}-f_2f_1^{\iota}),\\
 C=(f_3^{\iota}-f_3)/(f_1f_3^{\iota}-f_3f_1^{\iota}),& D=(f_1-f_1^{\iota})/(f_1f_3^{\iota}-f_3f_1^{\iota}),\\
 E_0=(f_3^{\iota}-f_3)/(f_2f_3^{\iota}-f_3f_2^{\iota}),& F_0=(f_2-f_2^{\iota})/(f_2f_3^{\iota}-f_3f_2^{\iota}).
\end{array}
\end{equation}
All numerators and denominators are nonzero by the choice of $U$, and
these functions are invariant under $\iota$ and satisfy
\begin{equation}\label{eq:ABCD}
 Af_1+Bf_2=1,\qquad Cf_1+Df_3=1,\qquad E_0f_2+F_0f_3=1.
\end{equation}
The norm functions $\operatorname{Nm}(f_1)=f_1f_1^{\iota}$,
$\operatorname{Nm}(f_2)=f_2f_2^{\iota}$ and
$\operatorname{Nm}(f_3)=f_3f_3^{\iota}$ are also invariant.
Thus they can be viewed as rational functions of the quotient
coordinate $X=x(Q-U)$ on $E/\iota\simeq\mathbb P^1$.

\subsubsection{Evaluation on the quotient}
For three nonzero rational functions $g_1,g_2,g_3$ of $X$, we write
\[
 g_1=c_1\prod_a(X-a)^{m_a},\qquad
 g_2=c_2\prod_b(X-b)^{n_b},\qquad
 g_3=c_3\prod_c(X-c)^{r_c}.
\]
Here the factors run over their finite zeros and poles. We put
\begin{equation}\label{eq:RP1}
 \mathcal R_{\mathbb P^1}(g_1,g_2,g_3)=
   \sum_{\substack{a,b,c\ \mathrm{pairwise\ distinct}}}
       m_a n_b r_c
       \left[\frac{c-a}{b-a}\right].
\end{equation}
\subsubsection{The formula on the elliptic curve}
We extend this rule linearly to the following sum of seven exterior
products of invariant functions:
\begin{align}
 V(f_1,f_2,f_3)={}&-\operatorname{Nm}(f_1)\wedge D\wedge B-C\wedge \operatorname{Nm}(f_3)\wedge B
             -2C\wedge D\wedge B\notag\\
            &+\operatorname{Nm}(f_2)\wedge F_0\wedge A+E_0\wedge \operatorname{Nm}(f_3)\wedge A
             +2E_0\wedge F_0\wedge A-A\wedge B\wedge \operatorname{Nm}(f_3).
\end{align}
Then we put
\begin{equation}
 \mathcal R_E(f_1,f_2,f_3)
   =S_E(Af_1,f_3)+S_E(Cf_1,B)-S_E(E_0f_2,A)+\mathcal R_{\mathbb P^1}\bigl(V(f_1,f_2,f_3)\bigr).
\end{equation}
Thus the arguments of the Bloch symbols in this expression are values
of specified rational functions and cross ratios of their zeros and
poles, and all of them are algebraic. By
Lemma~\ref{lem:RE}, this element of the
rational pre-Bloch group does not depend on the choice of the involution.

\begin{lemma}\label{lem:RE}
For nonzero rational functions $f_1,f_2,f_3$ on $E$,
\begin{equation}\label{eq:REh}
 \mathcal R_E(f_1,f_2,f_3)=h_E(f_1\wedge f_2\wedge f_3).
\end{equation}
\end{lemma}

\begin{proof}
If a factor is constant, then $h_E(f_1\wedge f_2\wedge f_3)=0$ by
\eqref{eq:hC}, and by multilinearity and alternation, the
same holds if a ratio of two factors is constant. In the remaining case,
by the degree formula in \eqref{eq:hC}, $h_E$ is invariant
under $\iota$, and its value on a wedge of invariant functions is twice
the corresponding value on $\mathbb P^1$. We first show that, in the third
exterior power,
\begin{align}\label{eq:wedge}
 f_1\wedge f_2\wedge f_3={}&
    Af_1\wedge Bf_2\wedge f_3
      +Cf_1\wedge Df_3\wedge B-E_0f_2\wedge F_0f_3\wedge A\notag\\
 &-f_1\wedge D\wedge B-C\wedge f_3\wedge B-C\wedge D\wedge B\notag\\
 &+f_2\wedge F_0\wedge A+E_0\wedge f_3\wedge A
       +E_0\wedge F_0\wedge A-A\wedge B\wedge f_3.
\end{align}
Expanding each product as a sum in the rationalized multiplicative group,
we have
\begin{align*}
 Af_1\wedge Bf_2\wedge f_3
 &=f_1\wedge f_2\wedge f_3+A\wedge f_2\wedge f_3
   +f_1\wedge B\wedge f_3+A\wedge B\wedge f_3,\\
 Cf_1\wedge Df_3\wedge B
 &=f_1\wedge f_3\wedge B+C\wedge f_3\wedge B
   +f_1\wedge D\wedge B+C\wedge D\wedge B,\\
 -E_0f_2\wedge F_0f_3\wedge A
 &=-f_2\wedge f_3\wedge A-E_0\wedge f_3\wedge A
   -f_2\wedge F_0\wedge A-E_0\wedge F_0\wedge A.
\end{align*}
The terms $f_1\wedge B\wedge f_3$ and $f_1\wedge f_3\wedge B$ cancel,
as do $A\wedge f_2\wedge f_3$ and $-f_2\wedge f_3\wedge A$, and the seven
correction terms cancel the other terms except $f_1\wedge f_2\wedge f_3$.
Thus we have \eqref{eq:wedge}.

We apply $h_E$ to this expansion. By
\eqref{eq:ABCD}, the first three terms
give $S_E(Af_1,f_3)+S_E(Cf_1,B)-S_E(E_0f_2,A)$. We average each remaining
term with its pullback by $\iota$; for instance, we have
\[
 2h_E(f_1\wedge D\wedge B)=h_E(\operatorname{Nm}(f_1)\wedge D\wedge B)
      =2h_{\mathbb P^1}(\operatorname{Nm}(f_1)\wedge D\wedge B).
\]
A wedge whose three factors are invariant has value twice its value
on $\mathbb P^1$. Thus these seven terms give precisely
$h_{\mathbb P^1}(V(f_1,f_2,f_3))$.
For pairwise distinct finite $a,b,c$, we put
$u=(X-a)/(b-a)$ and $v=(c-a)/(b-a)$. Then, after discarding constant
factors, we have
\[
 h_{\mathbb P^1}\bigl((X-a)\wedge(X-b)\wedge(X-c)\bigr)
   =h_{\mathbb P^1}\bigl(u\wedge(1-u)\wedge(u-v)\bigr)=[v]
\]
by \eqref{eq:hC}, since
$\operatorname{div}(u-v)=[v]-[\infty]$. As repeated factors have zero
wedge, factoring the three rational functions, we obtain exactly
\eqref{eq:RP1}. Thus we complete the proof of the lemma.
\end{proof}

\subsection{The class \texorpdfstring{$C_P$}{C\_P} in Bloch symbols}

We return to the elliptic curve $E/L$, the torsion point $P\in E[N]$,
and the CM endomorphism $\rho$ of Section~\ref{sec:eisenstein}. For the
finite sum below, we assume that $\bar\rho=-\rho$; this is the case used
in Section~\ref{sec:descent}, and it halves the number of
terms. Since the evaluation of $h_E$ below is multilinear and constant
factors do not contribute to it, we may replace the normalized functions
of the principal decomposition of Subsection~\ref{sec:construction} by
any functions with the same divisors. We use the functions $a_P$ and
$d_N$ below, whose divisors give $\psi_P$ and $\delta_N$ by
\eqref{eq:divdN} and \eqref{eq:psiP}. The graphs of $\rho$ and $\bar\rho$ pass through $E[N]^2$, where
the precycle defining $Z_{\psi_P}$ is not defined. Thus we translate them
by a torsion point outside $E[N]$; this does not change their rational
Chow classes (see the proof of
Proposition~\ref{prop:bP}).

We choose a short Weierstrass equation of $E$ over $L$ and define
\begin{equation}
 d_N=
 \begin{cases}
 \displaystyle\prod_{\{Q,-Q\}\subset E[N]\setminus\{0\}}
                      (x-x(Q)),&N\text{ odd},\\[3pt]
 \displaystyle y\prod_{\{Q,-Q\}\subset E[N]\setminus E[2]}
                      (x-x(Q)),&N\text{ even}.
 \end{cases}
\end{equation}
Here each product runs over the indicated unordered pairs $\{Q,-Q\}$.
From the divisors of $x-x(Q)$ and $y$, we have
\begin{equation}\label{eq:divdN}
 \operatorname{div}(d_N)
   =\sum_{Q\in E[N]}[Q]-N^2[0]=-N^2\delta_N.
\end{equation}
If $P\ne0$, we put $a_P=x-x(P)$. Then
$\operatorname{div}(a_P)=[P]+[-P]-2[0]$, including the case $P=-P$. Hence we have
\begin{equation}\label{eq:psiP}
 \psi_P=\frac12\operatorname{div}(a_P)
          -\frac{N^2}{N^4-1}\operatorname{div}(d_N).
\end{equation}
For $P=0$, only the second term is present.

Let $T_0\notin E[N]$ be a nonzero torsion point, for example, a point of
prime order not dividing $N$. For $\kappa\in\{\rho,\bar\rho\}$, we
consider the map
\[
 \theta_{\kappa,T_0}:E\longrightarrow E^2,
 \qquad Q\longmapsto(Q,\kappa Q+T_0).
\]
We put $\Gamma_{\kappa,T_0}=\theta_{\kappa,T_0}(E)$.
The function $d_N$ satisfies $d_N(-Q)=(-1)^{N-1}d_N(Q)$, and $a_P$
is even when $P\ne0$; this allows us to combine the sign changes in
the projector before applying $\mathcal R_E$.
For a nonzero rational function $p$ on $E$ and maps $u,v:E\to E$,
each the composite of an isogeny and a translation, such that $u-v$ and
$u+v$ are nonconstant, we put
\begin{equation}\label{eq:Hp}
\begin{aligned}
 H_p(u,v)
   ={}&\mathcal R_E\bigl(p(u),d_N(u-v),d_N(v)\bigr)\\
     &-\mathcal R_E\bigl(p(u),d_N(u+v),d_N(v)\bigr).
\end{aligned}
\end{equation}
Here $p(u)$ means $p\circ u$, and similarly for the other entries; since
$u$, $v$ and $u\pm v$ are surjective, all entries are nonzero rational
functions. We denote by $\mathrm{id}$ the identity map of $E$, and for
$S\in E[N]$, we put $v_S(Q)=\rho Q+T_0+S$. The linear parts of
$\mathrm{id}$, $v_S$ and $\mathrm{id}\pm v_S$ are $1$, $\rho$ and
$1\pm\rho$, which are nonzero since $\rho$ is nonreal, so the pairs
$(\mathrm{id},v_S)$ and $(v_S,\mathrm{id})$ satisfy the condition above.
For $P\ne0$, we define
\begin{equation}\label{eq:bP}
 b_P=\sum_{S\in E[N]}\left(
       -\frac{H_{a_P}(\mathrm{id},v_S)-H_{a_P}(v_S,\mathrm{id})}{4N^6}
       +\frac{H_{d_N}(\mathrm{id},v_S)}{N^4(N^4-1)}
       \right).
\end{equation}
For $P=0$, we omit the first term in parentheses. In
$H_{a_P}(v_S,\mathrm{id})$, the roles of the two maps are
interchanged, and its entries are $a_P(v_S)$, $d_N(v_S\mp\mathrm{id})$
and $d_N$.
Thus the formula uses at most $6N^2$ evaluations of $\mathcal R_E$
($2N^2$ when $P=0$).
Let $M\subset k=\overline L$ be a finite extension of $L$ containing $T_0$, the auxiliary
points $U$, and every argument of the symbols produced by
\eqref{eq:SC} and \eqref{eq:RP1} for
these finitely many triples.

\begin{proposition}\label{prop:bP}
Suppose that $\bar\rho=-\rho$, and let $Z_{\psi_P}$ be constructed
from a principal decomposition \eqref{eq:decomposition}
obtained by decomposing the three factors of
$\psi_P*\delta_N*\delta_N$ separately, using the normalized
functions $f_{A,B}$ of Subsection~\ref{sec:construction}. Then
the finite sum \eqref{eq:bP} belongs to
$\mathcal B_2(M)_{\mathbb Q}$ and satisfies
\begin{equation}\label{eq:PhibP}
 \Phi_M(b_P)=\operatorname{res}_{M/L}C_P
       \quad\text{in }\mathrm{CH}^2(M,3)_{\mathbb Q}.
\end{equation}
\end{proposition}

\begin{proof}
We first base-change the curve and the construction to $k=\overline L$,
and replace the graphs of $\rho$ and $\bar\rho$ by shifted graphs.
For $\kappa\in\{\rho,\bar\rho\}$, let $F_\kappa:E^2\to E$ be the map
$F_\kappa(Q_1,Q_2)=Q_2-\kappa Q_1$. Since the shifted graph
$\Gamma_{\kappa,T_0}$ is the fiber of $F_\kappa$ over $T_0$, we have
\[
 [\Gamma_{\kappa,T_0}]-[\Gamma_\kappa]
      =F_\kappa^*([T_0]-[0])=0
             \quad\text{in }\mathrm{CH}^1(E^2)_{\mathbb Q}.
\]
Here, if $T_0$ has order $m$, then $m([T_0]-[0])$ is principal on $E$,
and its pullback is zero in the rational Chow group. Moreover,
$\Gamma_{\kappa,T_0}$ is disjoint from $E[N]^2$: if $Q$ and
$\kappa Q+T_0$ were both $N$-torsion, then $T_0$ would be $N$-torsion.

With the decomposition in the proposition, the $e_\chi$-projections of the ordered graph symbols give the closed precycle on $E^2\setminus E[N]^2$ used to define
$Z_{\psi_P}$. Each ordered graph symbol in this precycle is the part over
$E^2\setminus E[N]^2$ of the graph of a morphism to $(\mathbb P^1)^3$,
intersected with the cube. Hence it meets
$\Gamma_{\kappa,T_0}\times\square^3$ properly, in the graph over
$\Gamma_{\kappa,T_0}$, and its pullback is the graph symbol of the
restricted functions. We show that the pullback of the precycle is admissible and has zero boundary. After applying $g\circ\theta_{\kappa,T_0}$, the
three maps to $E$ used to form the restricted triple have linear parts
among $\pm1$, $\pm\kappa$ and $\pm1\pm\kappa$, which are nonzero since
$\kappa$ is nonreal. Hence each map is an isogeny followed by a
translation, and the zero and pole sets of every restricted nonconstant
function are finite. This proves proper intersection with each face of codimension one. Each $f_{ij}$ has divisor supported on $E[N]$, and if
$(Q_1,Q_2)\in E^2$ and two of $-Q_1$, $Q_1-Q_2$ and $Q_2$ lie in $E[N]$,
then $(Q_1,Q_2)\in E[N]^2$. Since $G$ preserves $E[N]^2$, the pairwise
intersections of the zero and pole sets of each ordered triple occurring in the sum defining $e_\chi$ are contained in $E[N]^2$. The shifted graph is
disjoint from $E[N]^2$, so the restricted triples meet no face of
codimension two. Since the precycle has zero boundary and its pullback
meets all faces properly, the pullback has zero boundary as well, and
pulling back this precycle computes the pullback of $Z_{\psi_P}$. The rational equivalence above allows us to replace
$\Gamma_{\bar\rho}-\Gamma_\rho$ in \eqref{eq:CP} by
$\Gamma_{\bar\rho,T_0}-\Gamma_{\rho,T_0}$ without changing its product
with $Z_{\psi_P}$, since the product depends only on the rational
equivalence class of the divisor \cite{Bloch1986}.

Now we compute $\Phi_k^{-1}(\operatorname{res}_{k/L}C_P)$. The divisor
map is injective on $(k(E)^\times/k^\times)\otimes\mathbb Q$. For
$D=\psi_P$ or $\delta_N$, let $v_D$ be the unique element of this space
with divisor $D$. In the tensor products below, functions denote their
classes in this rationalized group. For $P\ne0$, by the decomposition and
\eqref{eq:divdN}--\eqref{eq:psiP}, we have
\begin{equation}\label{eq:tensor}
\begin{aligned}
 v_{\psi_P}\otimes v_{\delta_N}\otimes v_{\delta_N}
 &=\left(\frac{a_P}{2}-\frac{N^2d_N}{N^4-1}\right)
           \otimes\left(-\frac{d_N}{N^2}\right)
           \otimes\left(-\frac{d_N}{N^2}\right)\\
 &=\frac{a_P\otimes d_N\otimes d_N}{2N^4}
         -\frac{d_N\otimes d_N\otimes d_N}{N^2(N^4-1)}.
\end{aligned}
\end{equation}
For $P=0$, we omit the $a_P$ term in both right-hand sides.

Pullback along the three maps in an ordered triple, followed by the
exterior product and $h_E$, defines a trilinear map, whose value does not
change when an input function is multiplied by a nonzero constant. We
apply this map to \eqref{eq:tensor} and combine the
terms of the projector as follows.

Let $G_0\subset G$ be the subgroup of order eight generated by the two
inversions and the interchange. Each $g\in G$ can be written uniquely as
$h\circ t_{(R_1,R_2)}$ with $h\in G_0$ and $R_1,R_2\in E[N]$, where
$t_{(R_1,R_2)}$ is the translation by $(R_1,R_2)$. On the graph of
$\kappa$, the change of parameter $Q'=Q+R_1$ gives
\[
 t_{(R_1,R_2)}\theta_{\kappa,T_0}(Q)
   =(Q',\kappa Q'+T_0+R_2-\kappa R_1).
\]
Since translations have degree one, this change of parameter does not change the value of $h_E$, by the third identity in \eqref{eq:hC}.
Each $S=R_2-\kappa R_1\in E[N]$ occurs $N^2$ times, since for each $R_1$
there is exactly one $R_2=S+\kappa R_1$. Thus the sum over the
translations becomes $N^2$ times the sum over the shifts $T_0+S$.

We fix one such shift $t$, and put $u(Q)=Q$ and $v(Q)=\kappa Q+t$. For
$h\in G_0$, we put $(u_h,v_h)=h\circ(u,v)$. Since $a_P$ and $d_N$ are
even modulo constants, we may replace $p(-u_h)$ by $p(u_h)$. The four
elements of $G_0$ without interchange contribute $2H_p(u,v)$: the two with
equal signs each give the first evaluation in
\eqref{eq:Hp}, and the two with opposite signs
each give minus the second. Since the interchange has character $-1$, the
other four elements contribute $-2H_p(v,u)$. Thus, by
Lemma~\ref{lem:RE}, we have
\[
 \sum_{h\in G_0}\chi(h)\,
 h_E\bigl(p(-u_h)\wedge d_N(u_h-v_h)\wedge d_N(v_h)\bigr)
       =2\bigl(H_p(u,v)-H_p(v,u)\bigr).
\]
For $p=d_N$, interchanging the first and third wedge factors, we have
$H_{d_N}(v,u)=-H_{d_N}(u,v)$, and hence this sum is $4H_{d_N}(u,v)$.

It remains to combine the two graphs. Let $I_1(Q_1,Q_2)=(-Q_1,Q_2)$.
Since $\bar\rho=-\rho$, we have
\[
 I_1\circ \theta_{\rho,T_0}=\theta_{\bar\rho,T_0}\circ[-1].
\]
The projected precycle lies in the $\chi$-part, and $\chi(I_1)=-1$. Hence
its pullbacks to the two shifted graphs differ by a minus sign after
pullback by $[-1]$ on $E$, and this change of parameter does not change the value of $h_E$, by the third identity in \eqref{eq:hC}.
Thus subtracting the $\rho$-term from the $\bar\rho$-term multiplies the value of $h_E$ on the projected $\rho$-term by $-2$. Combining the translation
multiplicity, the sum over $G_0$, and the coefficients in
\eqref{eq:tensor} with the projector factor $1/(8N^4)$,
we obtain
\[
 \frac1{8N^4}\frac1{2N^4}N^2\cdot2\cdot(-2)=-\frac1{4N^6},
 \qquad
 \frac1{8N^4}\frac{-1}{N^2(N^4-1)}N^2\cdot4\cdot(-2)
       =\frac1{N^4(N^4-1)}.
\]
Thus the value of $h_E$ on the projected cycle is given by \eqref{eq:bP}, with $h_E$ in place of $\mathcal R_E$. Let $c_P$ be this cycle, formed from the normalized functions and pushed forward to $\operatorname{Spec}k$, with the difference $\Gamma_{\bar\rho,T_0}-\Gamma_{\rho,T_0}$ of shifted graphs. Then we have
$\operatorname{res}_{k/L}C_P=-[c_P]$ by \eqref{eq:CP}, and
$\Psi([c_P])=b_P$ by the evaluation above and \eqref{eq:REh}.
From \eqref{eq:PhiPsi}, we then have
\[
 \Phi_k^{-1}(\operatorname{res}_{k/L}C_P)
 =-\Psi(-[c_P])=\Psi([c_P])=b_P.
\]
Thus $b_P$ lies in the Bloch group over $k$, and
$\Phi_k(b_P)=\operatorname{res}_{k/L}C_P$.

Finally, we descend to $M$. Since $b_P$ is defined over $M$ and
$\bigwedge^2(M^\times\otimes\mathbb Q)$ injects into
$\bigwedge^2(k^\times\otimes\mathbb Q)$, the Bloch boundary of $b_P$,
which is zero over $k$, is zero over $M$. The equality of the two higher
Chow classes also descends to $M$. If a rational higher Chow
class becomes zero over $k$, then the finite cycle and boundary
witnessing this are defined over a finite extension $M'/M$.
Corestriction after restriction multiplies the class by $[M':M]$
\cite[Exercise~1.11 and Theorem~17.21]{MVW}, so the original rational
class is zero. Hence \eqref{eq:PhibP} follows from the
functoriality of $\Phi$. Thus we complete the proof of the proposition.
\end{proof}

\begin{lemma}\label{lem:isom}
Let $L$ be a number field and $N>3$. For $i=1,2$, let $E_i/L$ be an
elliptic curve with rational full $N$-torsion, let $P_i\in E_i[N]$, and
suppose that the nonreal CM endomorphisms $\rho_i,\bar\rho_i$ are defined
over $L$. Let $C_{P_i}$ be constructed with the principal decompositions of Subsection~\ref{sec:construction}. If an
origin-preserving isomorphism $f:E_1\to E_2$ over $L$ satisfies
$f(P_1)=P_2$, $f\rho_1=\rho_2f$, and $f\bar\rho_1=\bar\rho_2f$, then
$C_{P_1}=C_{P_2}$ in $\mathrm{CH}^2(L,3)_{\mathbb Q}$. Moreover, at fixed
$N$, the class $C_P$ does not depend on the basis of $E[N]$, the normalized functions, and the principal decompositions used in that construction.
\end{lemma}

\begin{proof}
We first work over $k=\overline L$. Let $T_1\notin E_1[N]$ be a torsion
point, and put $T_2=f(T_1)$. For corresponding endomorphisms
$\kappa_1,\kappa_2$, we have
$(f\times f)\theta_{\kappa_1,T_1}=\theta_{\kappa_2,T_2}f$ by the
compatibility with $f$. Since conjugation by $f\times f$ sends the
translation by $(R_1,R_2)$ to that by $(fR_1,fR_2)$ and preserves the
inversions and the interchange, it identifies the groups $G$ for $E_1$ and $E_2$, together with their characters $\chi$.

We denote by $\delta_{N,i}$ the divisor $\delta_N$ on $E_i$. Since $f$
sends $P_1$ to $P_2$ and maps $E_1[N]$ bijectively onto $E_2[N]$, it
pulls back $\psi_{P_2}$ and $\delta_{N,2}$ to $\psi_{P_1}$ and
$\delta_{N,1}$. Hence, by the injectivity of the divisor map modulo
constants, we have
\[
 f^*v_{\psi_{P_2}}=v_{\psi_{P_1}},\qquad
 f^*v_{\delta_{N,2}}=v_{\delta_{N,1}}.
\]
As in the proof of Proposition~\ref{prop:bP}, each projected cycle on the shifted graphs is evaluated by applying the ordered
pullbacks, the exterior product and $h_{E_i}$ to
$v_{\psi_{P_i}}\otimes v_{\delta_{N,i}}\otimes v_{\delta_{N,i}}$. By the
compatibility of the shifted graphs and the projectors, together with the
degree-one identity in \eqref{eq:hC}, these values are
equal. By \eqref{eq:CP} and Lemma~\ref{lem:Psi}, as in
the proof of Proposition~\ref{prop:bP}, each of them is
$\Phi_k^{-1}(\operatorname{res}_{k/L}C_{P_i})$. Applying $\Phi_k$, we see
that the restrictions of $C_{P_1}$ and $C_{P_2}$ to $k$ agree. Restriction
to $k$ is injective by the argument at the end of the proof of
Proposition~\ref{prop:bP}, so the equality holds over
$L$.

Finally, for $f=\mathrm{id}$, the same tensor shows that this value does not depend on the normalized functions and the decompositions. Since the group $G$ and its character depend only on
$E[N]$, changing the basis of $E[N]$ also leaves this value unchanged.
Thus $C_P$ does not depend on the stated choices at fixed $N$.
\end{proof}

\section{Descent to the ray class field}
\label{sec:descent}

So far, the arithmetic construction has been carried out over an
auxiliary number field, since the torsion points and the arguments of the
Bloch symbols are in general not defined over the ray class field $F$. In
this section, we descend to $F$ and prove Theorem~\ref{thm:main}. For
this, we use the reciprocity law of complex multiplication: a Galois
conjugate of our class is obtained by the same construction for another
ray class, and hence its regulator is again expressed by
$\Lambda_1(-1,\,\cdot\,)$ at that ray class
(Proposition~\ref{prop:conjugate}). In particular, all
conjugates by automorphisms fixing $F$ have the same regulator. Hence,
after we multiply the corresponding element of the Bloch group by the rational factor
$k_{K,\mathfrak f}$ of \eqref{eq:regCsigma}, the normalized transfer
to $F$, which averages these conjugates, gives the classes $\xi_\sigma$
with the required enhanced values.

\subsection{Complex multiplication reciprocity}

We now return to the imaginary quadratic field $K$, modulus
$\mathfrak f$, and ray class field $F=K_{\mathfrak f}$, with the complex
embedding and reciprocity convention fixed in the introduction. All
number fields below are subfields of $\mathbb C$ containing $K$,
$\overline K$ is the algebraic closure of $K$ in $\mathbb C$, and
regulators and enhanced values are taken at the inclusions into
$\mathbb C$ unless stated otherwise. To apply
the construction of Sections~\ref{sec:eisenstein}
and~\ref{sec:bloch}, we need a model of the pointed CM elliptic
curve over a number field. We first recall some standard facts about
such a model, mainly in order to fix conventions.

Let $H\subset F$ be the Hilbert class field. For nonzero fractional
ideals $\mathfrak b,\mathfrak c$ of $K$, let
$\sigma_{\mathfrak b}\in\operatorname{Gal}(H/K)$ be the reciprocity
automorphism associated with $\mathfrak b$.
With our convention, we have
$\sigma_{\mathfrak b}(j(\mathbb C/\mathfrak c))=j(\mathbb C/(\mathfrak b\mathfrak c))$
by \cite[Theorem~2.1]{Conrad}, and \cite[Theorem~2.3]{Conrad} describes
the Galois action on torsion points in terms of ideles.
For $\sigma\in\operatorname{Gal}(F/K)$, we put
$A=\operatorname{rec}^{-1}(\sigma)$, and take $\mathfrak a$ and $x_0$ as
in Section~\ref{sec:zeta}: $\mathfrak a$ is an integral
ideal in $A^{-1}$ prime to $\mathfrak f$, and $x_0\in\mathfrak a$
satisfies $x_0\equiv1\pmod{\mathfrak f}$. We consider the pointed complex
torus
\[
 E(\mathbb C)\simeq\mathbb C/(\mathfrak a\mathfrak f),
 \qquad P=x_0\bmod\mathfrak a\mathfrak f.
\]
Since $\Lambda=\mathfrak a\mathfrak f$ is an invertible $\mathcal O_K$-ideal,
the complex torus $\mathbb C/\Lambda$ has complex multiplication by the
maximal order $\mathcal O_K$. By the theory of complex multiplication, we
have $j(\mathbb C/\mathfrak a\mathfrak f)\in H$, and an elliptic curve
with this $j$-invariant admits a Weierstrass model over $H$;
see \cite[Theorem~2.1 and Corollary~2.2]{Conrad}. We choose such a model
$E/H$ and the analytic identification above. Note that the $j$-invariant
depends only on the ideal class of $\Lambda$, while the torsion point $P$
carries the additional information of the ray class.

Every endomorphism of $E$ is defined over $\overline H$, and we label the
endomorphisms by $\mu\in\mathcal O_K$ through their action
$\mu^*\omega=\mu\omega$ on a nonzero invariant differential $\omega$
defined over $H$. In characteristic zero, an endomorphism is determined by
this action. An automorphism $\widetilde{\gamma}$ fixing $K$ sends the
endomorphism labeled $\mu$ to the endomorphism of $E^{\widetilde{\gamma}}$
labeled $\widetilde{\gamma}(\mu)=\mu$, measured with
$\omega^{\widetilde{\gamma}}$. In particular, an automorphism fixing $H$
fixes every endomorphism of $E$, so the CM endomorphisms are defined over
$H$, and conjugation preserves the labels $\rho$ and $\bar\rho$.

We use $\rho=\sqrt{-d}$ with positive imaginary part and the scalar
$k_{K,\mathfrak f}$ from \eqref{eq:kKf}.
We choose an integer $N>3$ annihilating $P$ and an auxiliary finite
extension $L/F$ over which $E[N]$ is rational. With the decomposition of
Proposition~\ref{prop:bP}, the arithmetic construction
gives $C_\sigma=C_P\in \mathrm{CH}^2(L,3)_{\mathbb Q}$.

Combining Propositions~\ref{prop:regCP} and~\ref{prop:fourier}
with the normalization~\eqref{eq:kKf}, we obtain
in $\mathbb C/\mathbb Q(2)$,
\begin{equation}\label{eq:regCsigma}
\begin{aligned}
 R(k_{K,\mathfrak f}C_\sigma)
 &=\frac{k_{K,\mathfrak f}r}{2\pi i}G_\alpha(\tau)\\
 &=\frac{a_0}{w(2\pi i)}G_\alpha(\tau)\\
 &=2\pi i\Lambda_1(-1,\sigma).
\end{aligned}
\end{equation}

The construction is compatible with field conjugation. A function $f=f_{A,B}$ on $\mathcal E$, normalized at a torsion section $O'$ by $f(O')=1$, is sent
to $\widetilde{\gamma}(f)$, which has the conjugate divisor and satisfies
$\widetilde{\gamma}(f)(\widetilde{\gamma}O')=1$. The component $Y_N$ is
sent to the component attached to the conjugate value
of the Weil pairing. The rational coefficients of the decomposition are
unchanged, and conjugation identifies the groups $G$ and the characters $\chi$. Thus conjugation sends the closed precycle on
$\mathcal X\setminus\mathcal S$ to the precycle constructed from the
conjugate data. Since the extension in the localization sequence is
unique, and fiber pullback, graph products, and pushforward are compatible
with field isomorphisms, $\widetilde{\gamma} C_\sigma$ is the class
obtained by the same construction on
$(E^{\widetilde{\gamma}},P^{\widetilde{\gamma}})$. As we saw
above, conjugation by $\widetilde{\gamma}$ fixing $K$ preserves the labels $\rho$ and $\bar\rho$.

\begin{proposition}\label{prop:conjugate}
For every $\gamma\in\operatorname{Gal}(F/K)$ and every lift
$\widetilde{\gamma}\in\operatorname{Gal}(\overline K/K)$ of $\gamma$,
\begin{equation}\label{eq:conjugate}
 R(\widetilde{\gamma}^{-1}C_\sigma)=\frac{2\pi i}{k_{K,\mathfrak f}}\Lambda_1(-1,\sigma\gamma)
       \quad\text{in }\mathbb C/\mathbb Q(2).
\end{equation}
Here the regulator of $\widetilde{\gamma}^{-1}C_\sigma$ is taken at the fixed complex
embedding of its coefficient field $\widetilde{\gamma}^{-1}(L)$.
\end{proposition}

\begin{proof}
We first choose an idele representing $\widetilde{\gamma}$. Put
$B=\operatorname{rec}^{-1}(\gamma)$. Let $\mathfrak b$ be an integral
ideal in $B$ prime to $\mathfrak f$, and let $s_{\mathfrak b}$ be a finite
idele with associated ideal $\mathfrak b$ that is equal to one outside the
support of $\mathfrak b$. We also write $\operatorname{rec}$ for idelic
reciprocity with the same convention, that is, for the inverse of the
idelic Artin map of \cite{MilneCFT}; thus
$\operatorname{rec}(s_{\mathfrak b})|_F=\gamma$. Then the quotient
$\widetilde{\gamma}|_{K^{\mathrm{ab}}}\operatorname{rec}(s_{\mathfrak b})^{-1}$
fixes $F$. Idelic reciprocity is surjective onto
$\operatorname{Gal}(K^{\mathrm{ab}}/K)$, and the ideles whose image fixes
$F$ are the products of a principal idele and an element of the ray
subgroup $W_{\mathfrak f}$
\cite[Chapter~V, Propositions~4.6 and~4.7, Theorem~5.3, and Remarks~5.7(a) and~5.8]{MilneCFT}.
Thus the quotient is represented by the product of a principal idele, an
infinite component, and a finite unit idele $u=(u_v)$ with
$u_v\in\mathcal O_v^\times$ at every finite place and
$u_v\equiv1\pmod{\mathfrak f_v}$ for $v\mid\mathfrak f$. Here
$\mathcal O_v$ is the local integer ring, and a subscript $v$ on an ideal
denotes its extension to $\mathcal O_v$. The principal and complex
infinite factors have trivial image under $\operatorname{rec}$, so the
quotient is equal to $\operatorname{rec}(u)$. Put $s=s_{\mathfrak b}u$.
Then we have $(s)=\mathfrak b$ and
$\operatorname{rec}(s)=\widetilde{\gamma}|_{K^{\mathrm{ab}}}$; in the
usual normalization of the Artin symbol, which \cite{Conrad} uses, the
symbol of $s$ is therefore $\widetilde{\gamma}^{-1}|_{K^{\mathrm{ab}}}$.
We extend $\widetilde{\gamma}^{-1}$ to $\mathbb C$. Applying CM
reciprocity \cite[Theorem~2.3]{Conrad} to this inverse automorphism and
the idele $s$, we obtain the conjugate analytic lattice
$\mathfrak b^{-1}\mathfrak a\mathfrak f$, on which the conjugate point is $s^{-1}x_0$, where $s^{-1}$ acts on $K/\mathfrak a\mathfrak f$ through its finite components. Now we show that this point is represented by $x_0$ itself.
At $v\mid\mathfrak f$, we have $s_{\mathfrak b,v}=1$,
$s_v^{-1}-1\in\mathfrak f_v$, and $\mathfrak b_v=\mathcal O_v$. Since
$x_0\in\mathfrak a_v$, we have
\[
 (s_v^{-1}-1)x_0\in\mathfrak a_v\mathfrak f_v
                  =\mathfrak b_v^{-1}\mathfrak a_v\mathfrak f_v.
\]
At every other prime, we have $\mathfrak f_v=\mathcal O_v$, and both $x_0$
and $s_v^{-1}x_0$ lie in $\mathfrak b_v^{-1}\mathfrak a_v$; the first
inclusion follows from the integrality of $\mathfrak b$, and the second
from $(s_v)=\mathfrak b_v$. Thus $x_0-s^{-1}x_0$ lies in the new lattice
at every finite place.

Finally, we determine the ray class after rescaling. Since
$\mathfrak b+\mathfrak f=\mathcal O_K$, we can take
$0\ne c\in\mathfrak b$ with $c\equiv1\pmod{\mathfrak f}$. Put
\[
 \mathfrak a'=c\mathfrak b^{-1}\mathfrak a,
 \qquad x'_0=cx_0.
\]
Since $c\mathfrak b^{-1}\subset\mathcal O_K$ and $c$ is a unit at every
prime dividing $\mathfrak f$, the ideal $\mathfrak a'$ is integral and
prime to $\mathfrak f$. Moreover, we have $x'_0\in\mathfrak a'$ and
$x'_0\equiv1\pmod{\mathfrak f}$, and in the ray class group, we have
\[
 [\mathfrak a']=[\mathfrak b]^{-1}[\mathfrak a]
       =B^{-1}A^{-1}=(AB)^{-1},
\]
since $(c)$ has trivial ray class by $c\equiv1\pmod{\mathfrak f}$.
Multiplication by $c$ gives the pointed torus
$(\mathbb C/\mathfrak a'\mathfrak f,x'_0)$ and respects the CM action.
Hence the associated ray class is $AB$, whose image under
$\operatorname{rec}$ is $\sigma\gamma$.

We use the same integer $N$ for the conjugate pointed curve; it still annihilates the point. The resulting isomorphism of the pointed CM
curves, compatible with $\rho$ and $\bar\rho$, is defined over a finite
extension, and by Lemma~\ref{lem:isom}, the conjugate
class coincides with the class constructed from $(\mathfrak a',x'_0)$,
after restriction to a common coefficient field. Since
$k_{K,\mathfrak f}$ does not depend on the ray class and the scaling by
$c$, applying \eqref{eq:regCsigma} to $(\mathfrak a',x_0')$, we obtain
\eqref{eq:conjugate} for the given $\widetilde{\gamma}$. Thus we
complete the proof of the proposition.
\end{proof}

\subsection{Transfer and proof of the main theorem}

For each $\sigma\in\operatorname{Gal}(F/K)$, we form $b_P$ by
\eqref{eq:bP}, and take a finite extension $M/K$, normal
over $K$, that contains $L$ and all the points and values required before
Proposition~\ref{prop:bP}. Put
\[
 \beta_\sigma=k_{K,\mathfrak f}b_P=\sum_j c_j[z_j]
       \in\mathcal B_2(M)_{\mathbb Q},
\]
and define
\begin{equation}\label{eq:xi}
 \xi_\sigma=\frac1{[M:F]}\operatorname{cor}^{\mathcal B}_{M/F}(\beta_\sigma)
       \in\mathcal B_2(F)_{\mathbb Q}.
\end{equation}
Here $\operatorname{cor}^{\mathcal B}_{M/F}$ is the finite pushforward on
higher Chow groups transported through the isomorphisms $\Phi_M$ and
$\Phi_F$ of Corollary~\ref{cor:iso}.
Restriction is induced by the inclusion $F\hookrightarrow M$ and sends
each symbol $[z]$ to the same symbol over $M$. For
$H_0=\operatorname{Gal}(M/F)$, the two maps satisfy
\[
 \operatorname{cor}^{\mathcal B}_{M/F}\operatorname{res}_{M/F}
       =[M:F]\,\mathrm{id},\qquad
 \operatorname{res}_{M/F}\operatorname{cor}^{\mathcal B}_{M/F}
       =\sum_{h\in H_0}h.
\]
These identities follow from finite pushforward and flat pullback on
higher Chow groups \cite[Exercise~1.11 and Theorem~17.21]{MVW}.
By the first identity, rational restriction is injective, and by the
second, every invariant class descends by normalized transfer. Thus
$\mathcal B_2(F)_{\mathbb Q}$ is identified with
$\mathcal B_2(M)_{\mathbb Q}^{H_0}$, and $\operatorname{res}_{M/F}\xi_\sigma$
is the Galois average of $\beta_\sigma$.

\begin{proof}[Proof of Theorem~\ref{thm:main}]
We first compute the enhanced values of $\gamma^{-1}\xi_\sigma$. By
Proposition~\ref{prop:bP}, we have
\[
 \Phi_M(\beta_\sigma)=k_{K,\mathfrak f}\operatorname{res}_{M/L}C_\sigma.
\]
For $\gamma\in\operatorname{Gal}(F/K)$, let
$\widetilde{\gamma}\in\operatorname{Gal}(M/K)$ be a lift of $\gamma$. Since $M/K$ is
normal, we have $\widetilde{\gamma}^{-1}(L)\subset M$, and by the
naturality of $\Phi$, we have
\[
 \Phi_M(\widetilde{\gamma}^{-1}\beta_\sigma)
   =k_{K,\mathfrak f}\operatorname{res}_{M/\widetilde{\gamma}^{-1}(L)}
       (\widetilde{\gamma}^{-1}C_\sigma).
\]
Since restriction does not change the regulator at the fixed embedding,
applying Proposition~\ref{prop:conjugate} to an extension of
$\widetilde{\gamma}$ to $\overline K$ and using
Proposition~\ref{prop:regL2hat}, we have
\begin{equation}\label{eq:L2hat-beta}
 \sum_j c_j\widehat{\mathcal L}_2(\widetilde{\gamma}^{-1}z_j,\ell)
   =\widehat{\mathcal L}_2(\widetilde{\gamma}^{-1}\beta_\sigma)
   =2\pi i\Lambda_1(-1,\sigma\gamma)
       \quad\text{in }\mathbb C/\mathbb Q(2),
\end{equation}
where $\ell$ is a common $\mathbb Q$-linear logarithm lift on
$M^\times\otimes\mathbb Q$.

Since $H_0$ is normal in $\operatorname{Gal}(M/K)$, we have, by the
transfer identities above,
\[
 \operatorname{res}_{M/F}(\gamma^{-1}\xi_\sigma)
   =\frac1{|H_0|}\sum_{h\in H_0}h\widetilde\gamma^{-1}\beta_\sigma.
\]
Since restriction sends each symbol to the same symbol, it preserves
the enhanced value at the fixed embedding. Hence we have
\begin{align*}
 \widehat{\mathcal L}_2(\gamma^{-1}\xi_\sigma)
 &=\frac1{|H_0|}\sum_{h\in H_0}
         \widehat{\mathcal L}_2(h\widetilde{\gamma}^{-1}\beta_\sigma)\\
 &=2\pi i\Lambda_1(-1,\sigma\gamma).
\end{align*}
Here the last equality follows from \eqref{eq:L2hat-beta}
for the lifts $\widetilde{\gamma}h^{-1}\in\operatorname{Gal}(M/K)$ of
$\gamma$, whose inverses are $h\widetilde{\gamma}^{-1}$. Thus we have
\eqref{eq:main}.

Now we show that $\gamma^{-1}\xi_\sigma=\xi_{\sigma\gamma}$. By
\eqref{eq:main}, we have, for every
$\gamma'\in\operatorname{Gal}(F/K)$,
\[
 \widehat{\mathcal L}_2(\gamma'^{-1}\gamma^{-1}\xi_\sigma)
   =2\pi i\Lambda_1(-1,\sigma\gamma\gamma')
   =\widehat{\mathcal L}_2(\gamma'^{-1}\xi_{\sigma\gamma}).
\]
The embeddings of $F$ extending the fixed embedding of $K$ are precisely
the chosen embedding of $F$ composed with the automorphisms
$\gamma'^{-1}$, and the enhanced value of a class at such an embedding
is that of its $\gamma'^{-1}$-conjugate at the chosen embedding. Thus,
by Lemma~\ref{lem:inj}, we have
$\gamma^{-1}\xi_\sigma=\xi_{\sigma\gamma}$. In the same way, $\xi_\sigma$
does not depend on the auxiliary data, namely the ideal $\mathfrak a$ and
the element $x_0$, the model of $E$, the basis of $E[N]$, the level $N$,
the normalized functions and the principal decompositions, the shift
$T_0$, the points $U$, and the fields $L$ and $M$, since any two
constructions give the same enhanced values at every such embedding. Thus
we complete the proof of Theorem~\ref{thm:main}.
\end{proof}

\section*{Acknowledgements}
The author was supported by the National Science and Technology Council
(NSTC), Taiwan, under Grant No.~113-2115-M-002-007-MY3.

\end{document}